\documentclass[reqno]{amsart}

\usepackage[T1]{fontenc}
\usepackage[utf8]{inputenc}
\usepackage{lmodern}

\usepackage{amsmath,amssymb,amsthm,mathtools}
\usepackage{graphicx}
\usepackage{booktabs}
\usepackage{enumitem}
\usepackage{xcolor}
\usepackage{microtype}
\usepackage{nicefrac}
\usepackage{cite}
\usepackage[
  unicode,
  colorlinks=true,
  linkcolor=blue,
  citecolor=blue,
  urlcolor=blue
]{hyperref}
\usepackage{cleveref}
\usepackage[export]{adjustbox}

\usepackage{pifont}
\newenvironment{eqs} %
 { \begin{equation} \begin{aligned} } %
 { \end{aligned} \end{equation} \ignorespacesafterend } %

\theoremstyle{definition}
\theoremstyle{plain}
\makeatletter
\@ifclassloaded{beamer}{}{%
\newtheorem{theorem}{Theorem}[section]

\newtheorem{remark}[theorem]{Remark}
\crefname{assumption}{Assumption}{Assumptions}
}
\makeatother
\newtheorem*{remark*}{Remark}

\newtheorem*{acknowledgment*}{Acknowledgment}

\newcommand{\bpm}{\begin{pmatrix}}
\newcommand{\epm}{\end{pmatrix}}

\newcommand{\jmp}[1]{[\![ #1 ]\!]}

\newcommand{\avg}[1]{\{\!\!\{#1\}\!\!\}}
\newcommand{\dom}{\Omega} 

\newcommand{\Th}{{\calT_h}}
\newcommand{\Fh}{{\calF_h}}
\newcommand{\Fhi}{{\calF_h^i}}

\newcommand{\RT}{\calR\calT}

\renewcommand{\Re}{\operatorname{Re}}

\DeclareMathOperator{\Div}{div}

\DeclareMathOperator{\grad}{\nabla}

\newcommand{\inner}[1]{( #1 )}

\newcommand\restr[2]{{ \left.\kern-\nulldelimiterspace #1 \vphantom{\big|} \right|_{#2} }}

\newcommand{\IC}{\mathbb{C}}

\newcommand{\IN}{\mathbb{N}}

\newcommand{\IP}{\mathbb{P}}

\newcommand{\IR}{\mathbb{R}}

\newcommand{\IT}{\mathbb{T}}

\newcommand{\calF}{\mathcal{F}}

\newcommand{\calO}{\mathcal{O}}

\newcommand{\calR}{\mathcal{R}}

\newcommand{\calT}{\mathcal{T}}

\newcommand{\calX}{\mathcal{X}}

\usepackage[most]{tcolorbox}

\usepackage{tikzducks}
\makeatletter
\newlength{\maybeimage@width}
\newcommand{\maybeimage}[2][width=.8\linewidth]{%
  \IfFileExists{#2}{%
    \includegraphics[#1]{#2}%
  }{%
    \begingroup
    \setlength{\maybeimage@width}{.8\linewidth}%
    \let\Gin@ewidth\Gin@exclamation
    \setkeys{Gin}{#1}%
    \ifx\Gin@ewidth\Gin@exclamation\else
      \setlength{\maybeimage@width}{\Gin@ewidth}%
    \fi
      \fbox{%
        \parbox[c][4cm][c]{\dimexpr\maybeimage@width-2\fboxsep-2\fboxrule\relax}{%
          \centering
          \begin{tikzpicture}[scale=.6]
             \duck[glasses,tie=blue]
          \end{tikzpicture}

          Image not generated yet

          \scriptsize\nolinkurl{#2}
        }%
      }%
    \endgroup
  }%
}
\makeatother

\usepackage{pgfplots} 
\pgfplotsset{compat=newest}
\usepgfplotslibrary{colorbrewer,groupplots}
\pgfplotsset{
	discard if not/.style 2 args={
	x filter/.append code={
	\edef\tempa{\thisrow{#1}}
	\edef\tempb{#2}
	\ifx\tempa\tempb
	\else
	
	\fi
	}
	}
}
\usepackage{pgfplotstable} 
\usepackage{booktabs} 
\usepackage{colortbl}
\makeatletter
\pgfplotstableset{
	discard if not/.style 2 args={
	row predicate/.append code={
	\def\pgfplotstable@loc@TMPd{\pgfplotstablegetelem{##1}{#1}\of}
	\expandafter\pgfplotstable@loc@TMPd\pgfplotstablename
	\edef\tempa{\pgfplotsretval}
	\edef\tempb{#2}
	\ifx\tempa\tempb
	\else
	\pgfplotstableuserowfalse
	\fi
	}
	}
}
\makeatother

\pgfplotscreateplotcyclelist{paulcolors1}{%
	orange,every mark/.append style={solid,fill=orange},mark=square*,very thick,mark size=3pt\\%
	teal,densely dashed,every mark/.append style={solid,fill=teal},mark=*,very thick,mark size=3pt\\%
}
\pgfplotscreateplotcyclelist{paulcolors2}{%
	orange,every mark/.append style={solid,fill=orange},mark=square*,very thick,mark size=3pt\\%
	orange,densely dashed,every mark/.append style={solid,fill=orange},mark=*,very thick,mark size=3pt\\%
	teal,every mark/.append style={solid,fill=teal},mark=square*,very thick,mark size=3pt\\%
	teal,densely dashed,every mark/.append style={solid,fill=teal},mark=*,very thick,mark size=3pt\\%
    violet,every mark/.append style={solid,fill=violet},mark=square*,very thick,mark size=3pt\\%
    violet,densely dashed,every mark/.append style={solid,fill=violet},mark=*,very thick,mark size=3pt\\%
}

\pgfplotscreateplotcyclelist{paulcolors3}{%
	orange,every mark/.append style={solid,fill=orange},mark=square*,very thick,mark size=3pt\\%
	orange,densely dashed,every mark/.append style={solid,fill=orange},mark=*,very thick,mark size=3pt\\%
	orange,dotted,every mark/.append style={solid,fill=orange},mark=diamond*,very thick,mark size=3pt\\%
	teal,every mark/.append style={solid,fill=teal},mark=square*,very thick,mark size=3pt\\%
	teal,densely dashed,every mark/.append style={solid,fill=teal},mark=*,very thick,mark size=3pt\\%
	teal,dotted,every mark/.append style={solid,fill=teal},mark=diamond*,very thick,mark size=3pt\\%
	violet,every mark/.append style={solid,fill=violet},mark=square*,very thick,mark size=3pt\\%
	violet,densely dashed,every mark/.append style={solid,fill=violet},mark=*,very thick,mark size=3pt\\%
	violet,dotted,every mark/.append style={solid,fill=violet},mark=diamond*,very thick,mark size=3pt\\%
}

\pgfplotscreateplotcyclelist{paulcolors4}{%
	orange,every mark/.append style={solid,fill=orange},mark=square*,very thick,mark size=3pt\\%
	teal,every mark/.append style={solid,fill=teal},mark=*,very thick,mark size=3pt\\%
	violet,every mark/.append style={solid,fill=violet},mark=diamond*,very thick,mark size=3pt\\%
    cyan,every mark/.append style={solid,fill=cyan},mark=star,very thick,mark size=3pt\\%
}

\pgfplotscreateplotcyclelist{paulcolors8}{%
	orange,every mark/.append style={solid,fill=orange},mark=square*,very thick,mark size=3pt\\%
	teal,every mark/.append style={solid,fill=teal},mark=*,very thick,mark size=3pt\\%
	violet,every mark/.append style={solid,fill=violet},mark=diamond*,very thick,mark size=3pt\\%
    cyan,every mark/.append style={solid,fill=cyan},mark=star,very thick,mark size=3pt\\%
	orange,densely dashed,every mark/.append style={solid,fill=orange},mark=square*,very thick,mark size=3pt\\%
	teal,densely dashed,every mark/.append style={solid,fill=teal},mark=*,very thick,mark size=3pt\\%
	violet,densely dashed,every mark/.append style={solid,fill=violet},mark=diamond*,very thick,mark size=3pt\\%
    cyan,densely dashed,every mark/.append style={solid,fill=cyan},mark=star,very thick,mark size=3pt\\%
}

\pgfplotscreateplotcyclelist{paulcolors6}{%
	orange,every mark/.append style={solid,fill=orange},mark=square,very thick,mark size=3pt\\%
	teal,every mark/.append style={solid,fill=teal},mark=o,very thick,mark size=3pt\\%
	violet,every mark/.append style={solid,fill=violet},mark=diamond,very thick,mark size=3pt\\%
    cyan,every mark/.append style={solid,fill=cyan},mark=star,very thick,mark size=3pt\\%
    magenta,every mark/.append style={solid,fill=magenta},mark=triangle,very thick,mark size=3pt\\%
    black,every mark/.append style={solid,fill=black},mark=pentagon,very thick,mark size=3pt\\%
}

\title[DDM for non-conforming Helmholtz discretizations]{
    Comparing domain decomposition preconditioners for non-conforming Helmholtz discretizations
}

\author{Moritz Gallauner}
\address{
    Faculty of Mathematics\\
    University of Vienna\\
    Oskar-Morgenstern-Platz 1\\
    1090 Vienna, Austria
}
\email{a12315670@unet.univie.ac.at}

\author{Emile Parolin}
\address{
    Sorbonne Université, Université Paris Cité, CNRS, INRIA, Laboratoire
    Jacques-Louis Lions, LJLL, EPC ALPINES, 4 place Jussieu, Paris, F-75005,
    France
}
\email{emile.parolin@inria.fr}

\author{Paul Stocker}
\address{
    Faculty of Mathematics\\
    University of Vienna\\
    Oskar-Morgenstern-Platz 1\\
    1090 Vienna, Austria
}
\email{paul.stocker@univie.ac.at}

\author{Igor Voulis}
\address{
    Institut für Numerische und Angewandte Mathematik\\
    Georg-August-Universität Göttingen\\
    Lotzestr. 16-18\\
    37083 Göttingen, Germany
}
\email{voulis.sc@math.uni-goettingen.de}



\begin{document}

\begin{abstract}
    We compare additive and multiplicative domain decomposition preconditioners, without coarse correction, for three non-conforming polynomial discretizations of Helmholtz problems:
    discontinuous Galerkin, embedded Trefftz discontinuous Galerkin, and hybrid discontinuous Galerkin methods.
    The preconditioners are studied within stationary and Krylov iterative solvers.

    All three discretizations lead to complex-symmetric systems and local subproblems with inherited impedance-type boundary conditions.
    These properties enable the use of conjugate gradients and allow the local matrices to be obtained directly by restriction of the global matrix.
    Numerical experiments in two and three dimensions demonstrate promising performance for large-scale Helmholtz problems.
\end{abstract}

\maketitle

\section{Introduction}

Standard discretizations of time-harmonic wave propagation, as modeled by the Helmholtz equation, lead to indefinite and usually non-hermitian matrices, posing difficulties for iterative solvers~\cite{ErnstGander12,Gander2022}.
Additionally, the pollution effect, related to resolving high frequencies in the discretization, imposes mesh refinement or polynomial order increase and results in large linear systems~\cite{IB95,IB97,MS10,MS11}.
Efficient and robust preconditioners are therefore essential to solve the linear systems in practice.
Domain decomposition (DD) methods~\cite{Toselli2005,Dolean2015} offer an attractive approach, thanks in particular to their parallel nature.

A first class of DD methods for time-harmonic wave propagation problems are non-overlapping Optimized Schwarz Methods (OSM).
They were first introduced in~\cite{Despres1991}, see~\cite{Claeys2022c} for a short review, \cite{Pechstein2023} for a more comprehensive overview, and~\cite{GanderZhang19} for a review of sequential decompositions.
In the seminal work on the ultra-weak variational formulation (UWVF)~\cite{CD98}, diagonal element-wise blocks are inverted before the linear system is solved (see also~\cite{Pernet2023}).
This can be interpreted as an instance of OSM, and a convergence result of fixed-point iterative methods is moreover proved.
In the same spirit, the polynomial-based hybridized formulation of~\cite{ModaveChaumontFrelet23} employs Robin-type unknowns on facets and can be seen as a particular case of OSM, with a convergence result.
A numerical study on non-overlapping DD preconditioners (with two levels) for plane-wave Trefftz methods is available in~\cite{APZ15}, including multiplicative preconditioners.

The second class consists of overlapping methods, in particular the Optimised Restricted Additive Schwarz method (ORAS)~\cite{Graham2017,Gong2021,Gong2022b,Gong2023}.
Theoretical analysis of overlapping DD preconditioners for Helmholtz-type problems discretized with non-conforming methods is somewhat limited, in contrast to the works on elliptic problems~\cite{AA08,AAH06,AH11}.
Recently, two-level approaches for overlapping DD methods have been proposed, either with plane waves~\cite{Farhat2000,Farhat2005,APZ15}, or spectral coarse spaces~\cite{Hu2025,Ma2025b,Lu2025,Parolin2025,Galkowski2026,Dolean2026}.
Two-level additive and multiplicative overlapping DD preconditioners for plane wave DG systems are designed and compared in~\cite{HuLi16PWDG}.

A key feature in all DD methods for Helmholtz-type problems is the use of impedance (or Robin) boundary conditions in the subdomains.
Discontinuous Galerkin (DG) and more generally non-conforming methods are attractive in this setting because they expose the element interfaces where numerical transmission conditions can be imposed.
For Helmholtz problems this is more than a technical convenience: the flux terms can be read as local impedance couplings, so the discretization already contains the kind of Robin information that domain-decomposition preconditioners exploit.
A direct practical consequence of having built-in Robin-type boundary conditions in the formulation is that it is possible to simply restrict the global matrix to obtain local matrices for the local subproblems.
The well-posedness of the local problems, and the invertibility of the local matrices, is inherited from the global problem.
In contrast, for ORAS preconditioners used in conforming discretizations, the local matrices obtained by restriction are not necessarily invertible and therefore usually require assembly of separate, independent, local matrices.
This remark was the motivation for this paper.

We consider three different non-conforming discretizations of the Helmholtz equation to which we apply domain decomposition preconditioners.

The first one is a standard DG formulation based on the indefinite DG framework of \cite{MPS13}.
It uses element-wise polynomial approximation spaces coupled through numerical fluxes across the facets.
In the present comparison, it serves as the baseline method and as the underlying formulation from which the embedded Trefftz discretization is constructed.

The second one is an embedded Trefftz DG (TDG) formulation, as introduced in \cite{lozinski19,LS_IJMNE_2023}, and subsequently applied to a non-coercive DG formulation for Helmholtz in \cite{SV2026}.
It restricts the underlying polynomial DG space to an element-wise Trefftz-type subspace without requiring an explicit construction of Trefftz basis functions, thereby reducing the number of globally coupled unknowns.
These first two formulations are closely related to the Trefftz DG family for the Helmholtz equation, see~\cite{CD98,BM08,M11,HMP11,HMPS14,HMP16pwdg} and~\cite{HMP16} for a review.
An extensive review of the computational properties of Trefftz methods compared to other non-conforming methods has been done in \cite{LSZ_PAMM_2024}.

The third discretization considered is a Hybrid discontinuous Galerkin (HDG) formulation.
It introduces additional unknowns on the element facets and eliminates the element unknowns locally, resulting in a globally coupled system posed only on the mesh skeleton.
HDG formulations for Helmholtz problems are considered in~\cite{CLX13MLHDG,Leumuller2023,Leumuller2025,Huber2013,Monk2010}.
In the present comparison, HDG provides a computationally efficient reference point with a substantially different distribution of the globally coupled unknowns.

\textbf{Contributions.}
The main objective of this work is to compare the efficiency of several types of common preconditioners for these three \emph{non-conforming} polynomial-based discretizations.
The preconditioners considered here are additive and multiplicative versions based on overlapping and non-overlapping partitions.
We consider these preconditioners without any coarse correction and apply them within fixed-point and Krylov-based iterative solvers.
To the best of our knowledge, domain decomposition preconditioners have not previously been studied for an embedded Trefftz DG discretization.

The discretizations considered have two decisive advantages in this context.
First, all three formulations lead to complex-symmetric matrices, enabling the use of conjugate gradients as an iterative solver, which enjoys a short-recurrence relation compared to the GMRES solver usually used for Helmholtz problems.
Second, the local problems appearing in the preconditioners inherit impedance-type boundary conditions from the global formulations.
For the DG and TDG formulations, we compare both non-overlapping and overlapping partitions, while the HDG formulation naturally includes some overlap.

\textbf{Outline.}
The three non-conforming discretizations considered are described in \Cref{sec:discretizations}.
\Cref{sec:preconditioners} introduces the additive and multiplicative domain decomposition preconditioners studied.
Extensive numerical results in both 2D and 3D are then presented in \Cref{sec:numerics}. 
This includes a study of the parameters in the formulations and the discretization parameters on the efficiency of the preconditioners, strong and weak scaling tests, test cases with trapping effects, and a more challenging 3D test case on a realistic geometry.

\section{Non-conforming discretizations}
\label{sec:discretizations}

Let $\Omega\subset \IR^d$ be a bounded Lipschitz domain with outward unit
normal $n$ and three disconnected boundary components \(\Gamma_{D}\),
\(\Gamma_{N}\) and \(\Gamma_{R}\) such that
\(\partial\Omega = \Gamma_{D} \cup \Gamma_{N} \cup \Gamma_{R}\).
We consider the impedance boundary-value problem:
Find \(u\) such that
\begin{eqs}\label{eq:helmholtz}
    -\Delta u - \omega^2 u &= f &&\quad \text{in } \Omega,\\
    u &= 0 &&\quad \text{on } \Gamma_{D}, \\
    \grad u \cdot n &= 0 &&\quad \text{on } \Gamma_{N},\\
    \grad u \cdot n + i\omega u &= g &&\quad \text{on } \Gamma_{R},
\end{eqs}
for given \(f\), \(g\) and where $\omega>0$ denotes the wave number.
For every bounded Lipschitz domain $\dom$,
and for admissible data $f\in H^1(\dom)'$, 
and $g\in H^{-1/2}(\Gamma_{R})$,
standard well-posedness for the weak solution of \eqref{eq:helmholtz} is
guaranteed.

In the following we consider non-conforming discretizations of \eqref{eq:helmholtz}.
Let $\Th$ be a partition of $\Omega$ into non-overlapping Lipschitz elements $K$ with $\overline\Omega=\bigcup_{K\in\Th}\overline K$.
The set of all mesh faces is denoted by $\Fh$, with interior subset $\Fhi$ and boundary subsets $\calF_{h}^{D}$, $\calF_{h}^{N}$ and $\calF_{h}^{R}$ for \(\Gamma_{D}\), \(\Gamma_{N}\) and \(\Gamma_{R}\) respectively.
On an interior face $F=K^+\cap K^-$ we fix one normal orientation, written as $n$ and pointing from $K^+$ to $K^-$.
We use the broken spaces $\IP^p(\Th):=\prod_{K\in\Th}\IP^p(K)$, $\IP^p(\Fh):=\prod_{F\in\Fh}\IP^p(F)$, where $\IP^p$ denotes the space of polynomials of degree at most $p$.

For (complex-valued) arguments, we will denote the standard sesquilinear $L^2(D)$ product by $\inner{r,s}_D=\int_D r\bar s $, with complex conjugation on the second argument in the complex case.

\begin{table}[ht!]
\centering
\resizebox{\linewidth}{!}{%
\begin{filecontents*}{./numex/helmholtz_simplicial_counts.csv}
quantity,method,dim2_p1,dim2_p2,dim2_p3,dim2_p4,dim2_p5,dim2_p6,dim3_p1,dim3_p2,dim3_p3,dim3_p4,dim3_p5,dim3_p6
ndof,DG,3,6,10,15,21,28,4,10,20,35,56,84
ndof,TDG,3,5,7,9,11,13,4,9,16,25,36,49
ndof,HDG,3,6,9,12,15,18,4,12,24,40,60,84
nnze,DG,36,144,400,900,1764,3136,80,500,2000,6125,15680,35280
nnze,TDG,36,100,196,324,484,676,80,405,1280,3125,6480,12005
nnze,HDG,30,120,270,480,750,1080,56,504,2016,5600,12600,24696
\end{filecontents*}
\pgfplotstabletypeset[
    col sep=comma,
    string type,
    columns={
        quantity,method,
        dim2_p2,dim2_p3,dim2_p4,dim2_p5,dim2_p6,
        dim3_p2,dim3_p3,dim3_p4,dim3_p5,dim3_p6
    },
    columns/quantity/.style={
        column name={},
        column type=l,
        string replace={ndof}{$\mathrm{ndof}/N_{\mathrm{El}}$},
        string replace={nnze}{$\mathrm{nnze}/N_{\mathrm{El}}$},
    },
    columns/method/.style={column name={method}, column type=l},
    columns/dim2_p2/.style={column name={$p=2$}, column type=r},
    columns/dim2_p3/.style={column name={$p=3$}, column type=r},
    columns/dim2_p4/.style={column name={$p=4$}, column type=r},
    columns/dim2_p5/.style={column name={$p=5$}, column type=r},
    columns/dim2_p6/.style={column name={$p=6$}, column type=r},
    columns/dim3_p2/.style={column name={$p=2$}, column type=r},
    columns/dim3_p3/.style={column name={$p=3$}, column type=r},
    columns/dim3_p4/.style={column name={$p=4$}, column type=r},
    columns/dim3_p5/.style={column name={$p=5$}, column type=r},
    columns/dim3_p6/.style={column name={$p=6$}, column type=r},
    every head row/.style={before row={\toprule & & \multicolumn{5}{c}{2D simplices} & \multicolumn{5}{c}{3D simplices}\\ \cmidrule(lr){3-7}\cmidrule(lr){8-12}},after row=\midrule,},
    every row no 0/.style={before row=\rowcolor{gray!15},},
    every row no 2/.style={before row=\rowcolor{gray!15},after row=\addlinespace[0.6em],},
    every row no 3/.style={before row=\rowcolor{gray!15},},
    every row no 5/.style={before row=\rowcolor{gray!15},},
    every last row/.style={after row=\bottomrule},
]{./numex/helmholtz_simplicial_counts.csv}}
\vspace{.5em}
\caption{
    Total number of degrees of freedom (\(\mathrm{ndof}\)) and total number of non-zero entries (\(\mathrm{nnze}\)) divided by the number of mesh elements (\(N_{\mathrm{El}}\)) (assuming periodic boundary conditions to remove boundary effects).
    For HDG only the globally coupled condensed trace system is counted; the discontinuous highest-order trace modes are eliminated locally.
}
\vspace{-1em}
\label{tab:helmholtz-simplicial-counts}
\end{table}

In \Cref{tab:helmholtz-simplicial-counts} we compare the relevant computational quantities of the methods discussed below. 
We consider the degrees of freedom ($\mathrm{ndof}$) and the number of non-zero entries of the global linear system ($\mathrm{nnze}$) per number of elements (\(N_{\mathrm{El}}\)).   
For HDG we display the variables for the condensed system only.

\subsection{DG discretization}
\label{sec:dg}

We use the formulation of \cite{MPS13}, which may be viewed as a symmetric interior-penalty-type DG formulation with an added stabilization of the jump of the normal derivative.
The parameter $\delta$ blends two consistent treatments of the impedance boundary condition: for $\delta=0$, one recovers its standard weak Robin formulation, whereas $\delta=1$ yields an alternative residual-based boundary formulation.
When restricted to elementwise Trefftz spaces, the formulation reduces, after integration by parts, to the skeleton formulation studied in \cite{M11,HMP11}.
This family of plane-wave DG formulations generalizes the ultra-weak variational formulation introduced in \cite{CD98}.

For a piecewise smooth function $v$ we write $v^\pm:=v|_{K^\pm}$ and set
\begin{equation*}
    \avg{v} := \frac{1}{2}(v^+ + v^-), \quad \jmp{v} := v^+ - v^-.
\end{equation*}
For the DG discretization, we consider the discrete space $V_h=\IP^p(\Th)$.
The DG formulation reads: Find $u_h\in V_h$ such that
\begin{equation}\label{eq:dg}
    a_h(u_h,v_h) = \ell_h(v_h) \quad \forall v_h\in V_h,
\end{equation}
with the sesquilinear form $a_h(\cdot,\cdot)$ given by
\begin{eqs}\label{eq:ah}
    a_h(u,v) :=\ &\inner{\grad u,\grad v}_\Th - \inner{\omega^2 u, v}_\Th  
    - \inner{\avg{\grad u\cdot n},\jmp{v}}_\Fhi 
    - \inner{\jmp{u},\avg{\grad v\cdot n}}_\Fhi 
    \\
    &+ i \omega\alpha \inner{\jmp{u},\jmp{v}}_\Fhi 
    - \frac{\beta}{i\omega} \inner{\jmp{\grad u \cdot n},\jmp{\grad v \cdot n}}_\Fhi
    \\
    &
    - \delta \inner{u , \grad v\cdot n}_{\calF^R_h}
    - \delta \inner{\grad u\cdot n , v}_{\calF^R_h}
    - \frac{\delta}{i\omega} \inner{\grad u \cdot n, \grad v \cdot n}_{\calF^R_h}
    \\
    &+ i\omega (1-\delta) \inner{u, v}_{\calF^R_h}
    - \inner{\nabla u \cdot n,v}_{\calF^D_h} - \inner{u,\nabla v \cdot n}_{\calF^D_h}
    + \inner{\alpha u,v}_{\calF^D_h}
\end{eqs}
and $\ell_h(\cdot)$ given by
\begin{eqs}\label{eq:lh}
    \ell_h(v_h) :=&\ \inner{f, v_h}_\Th
    - \frac{\delta}{i\omega} \inner{g, \grad v_h \cdot n}_{\calF^{R}_{h}}.
\end{eqs}

\begin{remark}[Choice of parameters]
    On the choice of the parameters $\alpha,\beta,\delta>0$ we refer to \cite{MPS13}
    and to Sec.~\ref{sec:param-disc}.
    The above formulation coincides with the Ultra Weak Variational Formulation
    of~\cite{CD98} (which employs plane waves) if, and only if, \(4\alpha\beta
    = 1\) and \(\delta = 1/2\).
\end{remark}

\begin{remark}[Complex-symmetric linear system]
Although the variational products are written sesquilinearly, the real local polynomial bases used in the implementation make the element matrices transpose-symmetric.
The only complex quantities in them are scalar coefficients such as $i\omega\alpha$ and $-\beta/(i\omega)$, so the assembled matrix is complex symmetric ($A^T=A$), but in general not Hermitian.
This allows the use of a complex-symmetric direct solver, using an $LDL^T$ factorization.
Additionally, for symmetric preconditioners, one can employ a complex-symmetric variant of the conjugate-gradient method as iterative solver.
\end{remark}

\subsection{Embedded Trefftz DG discretization}
\label{sec:tdg}

Trefftz DG methods reduce the number of DG unknowns by restricting the discrete space to (approximate) solutions to the underlying PDE.
For the Helmholtz equation, a basis of exact solutions can be constructed from plane waves~\cite{HMP16}.
However, the non-polynomial basis functions can lead to ill-conditioned systems, especially for high-order approximations.
Additionally, standard Trefftz methods encounter difficulties when dealing with inhomogeneous equations or varying coefficients.

One way to overcome these issues and generalize the framework is to weaken the Trefftz condition. 
The embedded Trefftz construction was introduced in \cite{lozinski19,LS_IJMNE_2023} as a general approach to obtain Trefftz-like approximation spaces from standard polynomial finite element spaces without explicitly constructing problem-dependent Trefftz basis functions.
It has been applied to a wide range of wave and fluid dynamics problems (even non-linear) in \cite{LLS_NM_2024,SV2026, PIMPS_ARXIV_2026, GPSV_ARXIV_2026, SVLL1_ARXIV_2026, SVLL2_ARXIV_2026} and a framework for its analysis is available in \cite{LLSV_ARXIV_2024}.

In the case of the Helmholtz equation we will consider the local space
\begin{equation*}
\IT^p(K) := \{v_h \in \IP^p(K) \ | \ \inner{(-\Delta - \omega^2)v_h, q_h}_K = 0, \forall q_h \in \IP^{p-2}(K)\}.
\end{equation*} 
We do not compute the basis functions of this space explicitly. 
The kernel is computed by applying a singular value decomposition. 
An (orthogonal) basis of the local Trefftz space is then given by the corresponding vectors of the basis transformation matrix.
The method and concise construction is detailed in \cite{LS_IJMNE_2023}.

The Trefftz-DG formulation then simply swaps the discrete space $V_h$ in \eqref{eq:dg} with the Trefftz space 
$\IT^p := \Pi_{K \in \Th} \IT^p(K)$, and reads: Find $u_h \in \IT^p$ such that
\begin{equation}\label{eq:tdg}
    a_h(u_h,v_\IT) = \ell_h(v_\IT) \quad \forall v_\IT \in \IT^p,
\end{equation}
with $a_h(\cdot,\cdot)$ and $\ell_h(\cdot)$ as in \eqref{eq:ah} and \eqref{eq:lh}, respectively.

The embedded Trefftz DG method can also treat inhomogeneous equations by computing a particular solution in the polynomial space, and then solving for the homogeneous part in the Trefftz space.

\subsection{Hybrid DG discretization}
\label{sec:hdg}

HDG discretizations exploit the concept of static condensation to drastically reduce the system complexity, see~\cite{Monk2010,Huber2013,Leumuller2023,Leumuller2025}.
Element coupling is shifted to newly introduced facet spaces, allowing to condense all volume degrees of freedom.
Hence, only a smaller linear system involving the skeleton degrees of freedom has to be solved.
Such a method is obtained by discretizing the first-order mixed formulation:
Find \((u,\sigma)\) such that:
\begin{eqs}\label{eq:mixed_form}
i \omega \sigma - \nabla u &= 0 \quad &&\text{in } \Omega, \\
i\omega(-\Div \sigma + i \omega u) &= f \quad &&\text{in } \Omega, \\
u &= 0 \quad && \text{on } \Gamma_{D}, \\
\sigma \cdot n &= 0 \quad && \text{on } \Gamma_{N},\\
i\omega(\sigma \cdot n + u) &= g \quad &&\text{on } \Gamma_{R},
\end{eqs}
which is equivalent to \eqref{eq:helmholtz} by introducing the flux variable $\sigma := -i \omega^{-1}\nabla u$.

We define the composite discrete space 
\begin{equation}\label{eq:Xh}
\calX_h := \IP^p(\Th) \times \RT^p(\Th) \times \IP^p_D(\Fh) \times \IP^p(\Fh).
\end{equation}
Here, $\RT^{p}(\Th) := \prod_{K\in\Th}\RT^{p}(K),$ denotes the broken Raviart-Thomas space of order $p$.
And $\IP^p_D(\Fh)$ denotes the facet space of polynomials that vanish on $\Gamma_D$.
Its local degrees of freedom are given by moments of the normal component on each facet, 
supplemented by suitable interior moments, no continuity of the normal traces is imposed across facets.

The HDG formulation for \eqref{eq:mixed_form} reads: Find $(u,\sigma,\hat{u},\hat{\sigma}) \in \calX_h$ such that 
\begin{equation}\label{eq:hdg}
    b_h(u,\sigma,\hat{u},\hat{\sigma};v,\tau,\hat{v},\hat{\tau}) = i\omega^{-1} \inner{f,v}_\Omega + i\omega^{-1} \inner{g,\hat{v}}_{\Gamma_{R}},
\end{equation}
for all $(v,\tau,\hat{v},\hat{\tau}) \in \calX_h$, where $b_h$ is the sesquilinear form 
\begin{eqs}\label{eq:bh}
b_h(u,\sigma,\hat{u},\hat{\sigma};v,\tau,\hat{v},\hat{\tau}) := 
\sum_{K \in \Th} &\bigg(i\omega \inner{\sigma,\tau}_K + \inner{u,\Div \tau}_K +\inner{\Div \sigma,v}_K \\
    &- i\omega \inner{u,v}_K -\inner{\hat{u},\tau\cdot n_{K}}_{\partial K} - \inner{\sigma \cdot n_{K}, \hat{v}}_{\partial K} \\
&-\alpha\inner{[u],[v]}_{\partial K} + \beta \inner{[\sigma]_n,[\tau]_n}_{\partial K}\bigg) \\
&-\inner{\hat{u},\hat{v}}_{\Gamma_{R}}.
\end{eqs}
Here, \(n_{K}\) is the outward normal vector to \(K\) and $\alpha, \beta$ are positive stabilization parameters and the jumps with respect to the facet variables are given by
\begin{equation*}\label{eq:hdg_jumps}
    [u] := u - \hat{u}, \quad [\tau]_n := \tau \cdot n_{K} - (n \cdot n_{K}) \hat{\tau}.
\end{equation*}

\begin{remark}[Stabilization parameters]
The HDG formulation \eqref{eq:hdg} allows for an $h$-independent choice of stabilization parameters, 
$\alpha = \calO(1)$ and $\beta = \calO(1)$. 
Stabilizing in this way is beneficial for preconditioned iterative solvers, 
at the cost of complicating the analysis \cite{Leumuller2023}. 
\end{remark}

\begin{remark}[Quality of approximation]
We use projected jumps to reconstruct the primal variable to one degree higher approximation, for details we refer to \cite{Leumuller2023,Huber2013,L_MTH_2010}.
Hence the dual variable, and their facet trace, can be chosen one degree lower.
In the implementation the highest-order modes of $\hat{u}$ are discontinuous between neighboring elements and are also eliminated locally. 
Thus, the globally coupled condensed system contains two $\IP^{p-1}(F)$ trace fields per facet \(F\).
The resulting approximation properties match those of \eqref{eq:dg} with order $p$, 
while the condensed global system only retains facet spaces of order $p-1$.
\end{remark}

\section{Domain decomposition preconditioners}
\label{sec:preconditioners}

All formulations above yield a square linear system
\begin{equation}
    Ax = b,
\end{equation}
to be solved for \(x \in \IC^{N}\),
with matrix \(A \in \IC^{N \times N}\)
and right-hand side \(b \in \IC^{N}\),
where \(N\) is the total number of degrees of freedom.
We next consider various domain decomposition preconditioners to obtain efficient solvers in iterative methods.

\subsection{Additive Domain Decomposition Method}

\subsubsection{Non-overlapping method (block Jacobi)}

Block Jacobi uses a non-overlapping partition of the unknowns into blocks (e.g.\ elementwise blocks or subdomain patches) and builds a block-diagonal approximation of the system matrix.
In each preconditioning step, the residual is restricted to each block and a local problem is solved independently (typically by a direct solver on the block).
The global correction is then formed by assembling the sum of these local block corrections. 

More precisely, assume that we have a non-overlapping partition of the degrees of freedom into $J$ blocks, where the $j$-th block is of dimension $N_j$ such that the total number of degrees of freedom is $N = \sum_{j} N_j$.
Let the rectangular matrices $E_j \in \{0,1\}^{N \times N_j}$, be embedding matrices, mapping local to global degrees of freedom such that $E_j^T E_i = \mathbf 0\ \forall j\neq i$ (non-overlap condition).
The diagonal blocks are then defined by 
$A_j = E_j^T A E_j$,
and the preconditioner $M \in \IC^{N \times N}$ is given by $M = \text{diag}(A_1, \dots, A_J)$,
or equivalently
\begin{equation}\label{eq:additive}
M^{-1} = \sum_{j=1}^J E_j A_j^{-1} E_j^T.
\end{equation}
Note that the \emph{additive} structure allows computations to be parallelized easily.

\textbf{DG/TDG formulation.} 
The structure of DG/TDG lends itself naturally to the non-overlapping preconditioner. 
Consider a non-overlapping partition of the domain such that $\overline{\Omega}=\cup \overline{\Omega_j}$. 
Then the definition above corresponds to the blocks in the preconditioner being the degrees of freedom associated
to the elements contained in $\Omega_j$. 

\textbf{HDG formulation.}
With the HDG formulation, there is no natural way of assigning the degrees of freedom from the facet unknowns $(\hat u, \hat \sigma)$ at subdomain interfaces to a particular subdomain. 
Hence no non-overlapping preconditioners will be considered for HDG.

\subsubsection{Overlapping methods (Additive Schwarz)}

The preconditioner introduced above can also be defined and used with overlapping partitions of the degrees of freedom, often subordinate to a geometric overlapping partition of the domain.
The overlapping subdomains are denoted by \(\widetilde{\Omega}_{j}\).
In this case, \(E_{j}^{T}E_{i} \neq \mathbf{0}\) for subdomains \(\widetilde{\Omega}_{j}\) and \(\widetilde{\Omega}_{i}\) sharing some degrees of freedom.
The definition of the additive preconditioner in~\eqref{eq:additive} remains valid for overlapping partitions, but it is no longer block diagonal.

\textbf{HDG formulation.}
After condensation the global degrees of freedom of the method are only associated to the skeleton $\Fh$. 
Consider the same non-overlapping partition $\overline{\Omega}=\cup \overline{\Omega_j}$. 
The degrees of freedom in the \(j\)-th block are associated to some $F\in \Fh$ such that $F\subset\partial K$ for $K \in\Omega_j$. 
The overlap then naturally corresponds to the degrees of freedom associated with $F \in \Fh$ such that $F \subset \partial\Omega_j$.

\textbf{DG/TDG formulation.} 
The minimal overlap between subdomains \(\Omega_{j}\) and \(\Omega_{i}\) consists in the elements that contain a shared facet \(F \subset \overline{\Omega_{j}} \cap \overline{\Omega_{i}}\).
This type of overlap is chosen to mimic the behavior of HDG.

\subsubsection{Weighting (Restricted Additive Schwarz)}

The preconditioner in the overlapping case can be considered to over-correct in the overlap regions, which is the motivation to introduce a weighted version.
More precisely, we assume to have partition-of-unity (PoU) diagonal matrices \(D_{j} \in \IR^{N_{j} \times N_{j}}\) for each subdomain \(\widetilde{\Omega}_{j}\) such that 
\begin{equation}
    I = \sum_{j=1}^{J} E_{j} D_{j} E_{j}^{T}.
\end{equation}
Then the preconditioner in~\eqref{eq:additive} is modified as follows
\begin{equation}\label{eq:preaddpou}
M^{-1} = \sum_{j=1}^J E_j D_{j} A_j^{-1} E_j^T.
\end{equation}

\textbf{DG/TDG formulation.} 
Degrees of freedom in an element $K$ in the overlap are weighted in the diagonal entry of $D_{j}$ by $1/\#\{\widetilde{\Omega}_j : K\in\widetilde{\Omega}_j\}$, i.e.\ the inverse of the number of overlapping domains containing the element.

\textbf{HDG formulation.}
Degrees of freedom on shared facets between subdomains are weighted by $1/2$. 

\subsection{Multiplicative Domain Decomposition Method}

\subsubsection{Non-overlapping method (alternating block Gauss--Seidel)}

Block Gauss--Seidel is the block-triangular counterpart of block Jacobi: the same non-overlapping blocks are used, but the block solves are applied sequentially, updating the approximation after each block correction.
Alternating with the lower and upper triangular parts for symmetry, the preconditioner corresponds to a forward then backward sweep through the blocks, and the order can affect convergence.

Decomposing $A$ into $L+D+U$, where $L$ denotes the block-lower-triangular part, 
$D$ denotes the block-diagonal and $U$ denotes the block-upper-triangular part, 
one can write the preconditioner as
\begin{equation*}
    M^{-1} = (U+D)^{-1}(L+D)^{-1}.
\end{equation*}
This block-lower-triangular (resp. block-upper-triangular) system is solved by using (blockwise) forward (resp. backward) substitution. 
The preconditioner can be expressed in a more general way by using the embedding matrices $E_i$,  
\begin{equation}\label{eq:premult}
    M^{-1} = M^{-1}_{B} M^{-1}_{F},
    \ \text{where}\ 
    \begin{aligned}
        & M^{-1}_{F} =
        \sum_{j=1}^{J} E_j A_j^{-1} E_j^T \prod_{k=1}^{j-1} (I - A E_{j-k} A_{j-k}^{-1} E_{j-k}^T),
        \\
        & M^{-1}_{B} =
        \sum_{j=0}^{J-1} E_{J-j} A_{J-j}^{-1} E_{J-j}^T\!\! \prod_{k=J-j+1}^{J}\!\! (I - A E_{k} A_{k}^{-1} E_{k}^T),
    \end{aligned}
\end{equation}
where local solves can again be computed directly. 
Compared to block Jacobi, stationary iterative methods usually converges faster thanks to the immediate usage of updated information. 
However, the \emph{multiplicative} nature of this preconditioner allows for significantly less parallelism than additive methods.

\subsubsection{Overlapping methods (Multiplicative Schwarz and Restricted Multiplicative Schwarz)}

Overlapping partitions with or without weighting are also considered for these multiplicative preconditioners. 
We use the same overlap and partition of unity matrices as for the additive versions.
The partition of unity matrices $D_j$ are inserted into the preconditioner~\eqref{eq:premult} as factor of $A_j^{-1}$ in the same way as in~\eqref{eq:preaddpou}.

\subsection{Remarks on the local solves}

The local matrices \(A_{j}\) in all above preconditioners are obtained by restriction of the global matrix \(A\).
In conforming methods, the local solves then correspond to imposing non-homogeneous Dirichlet boundary conditions on subdomain interfaces.
On the contrary, in the considered non-conforming formulations, the local solves incorporate non-homogeneous Robin-type boundary conditions.
This is essential since simple stationary (fixed-point) iterative methods typically fail to converge with Dirichlet boundary conditions on subdomain interfaces for Helmholtz-type problems (even with overlap), see~\cite{Despres1991} and~\cite[Sec.~2.5]{Dolean2015}.
Additionally, in the non-overlapping case, Robin boundary conditions are necessary for convergence even for symmetric positive definite problems~\cite{Lions1990}.
As we shall see in~\Cref{sec:param-disc}, it is then remarkable that all non-conforming methods considered here can converge in stationary methods without modifying the local matrices thanks to Robin-type boundary conditions naturally present in these formulations, as already observed in~\cite{APZ15}.

To motivate this fact, we first note that the DG sesquilinear form \eqref{eq:ah} can be interpreted as a weighted Robin coupling of the interior jumps.
To see this, define the weighted Robin jumps
\[
    \calR_{\pm}^{\alpha,\beta} w
    :=
    \sqrt{\beta}\,\jmp{\partial_{n} w}
    \pm i\omega\sqrt{\alpha}\,\jmp{w}.
\]
Then the jump-jump part on $F\in\Fh$ satisfies the identity
\begin{equation*}\label{eq:robin-jump-penalty}
    i\omega\alpha
    \inner{ \jmp{u}, \jmp{v} }_F 
    -\frac{\beta}{i\omega}
    \inner{ \jmp{\partial_{n}u}, \jmp{\partial_{nF}v} }_F
    \!=\!
    \frac{i}{2\omega}\!\! \left(
          \inner{ \calR_+^{\alpha,\beta}u, \calR_+^{\alpha,\beta}v }_F
        + \inner{ \calR_-^{\alpha,\beta}u, \calR_-^{\alpha,\beta}v }_F
        \right). 
\end{equation*}
Thus the Dirichlet and Neumann jump contributions in \eqref{eq:ah} combine into a coupling of the two weighted Robin mismatch components, rather than acting as two unrelated penalties.
The same observation holds true for the jump-jump terms for the first-order formulation in HDG.

\begin{figure}[ht!]
\centering
\maybeimage[width=.32\linewidth, trim={250 80 180 80}, clip]{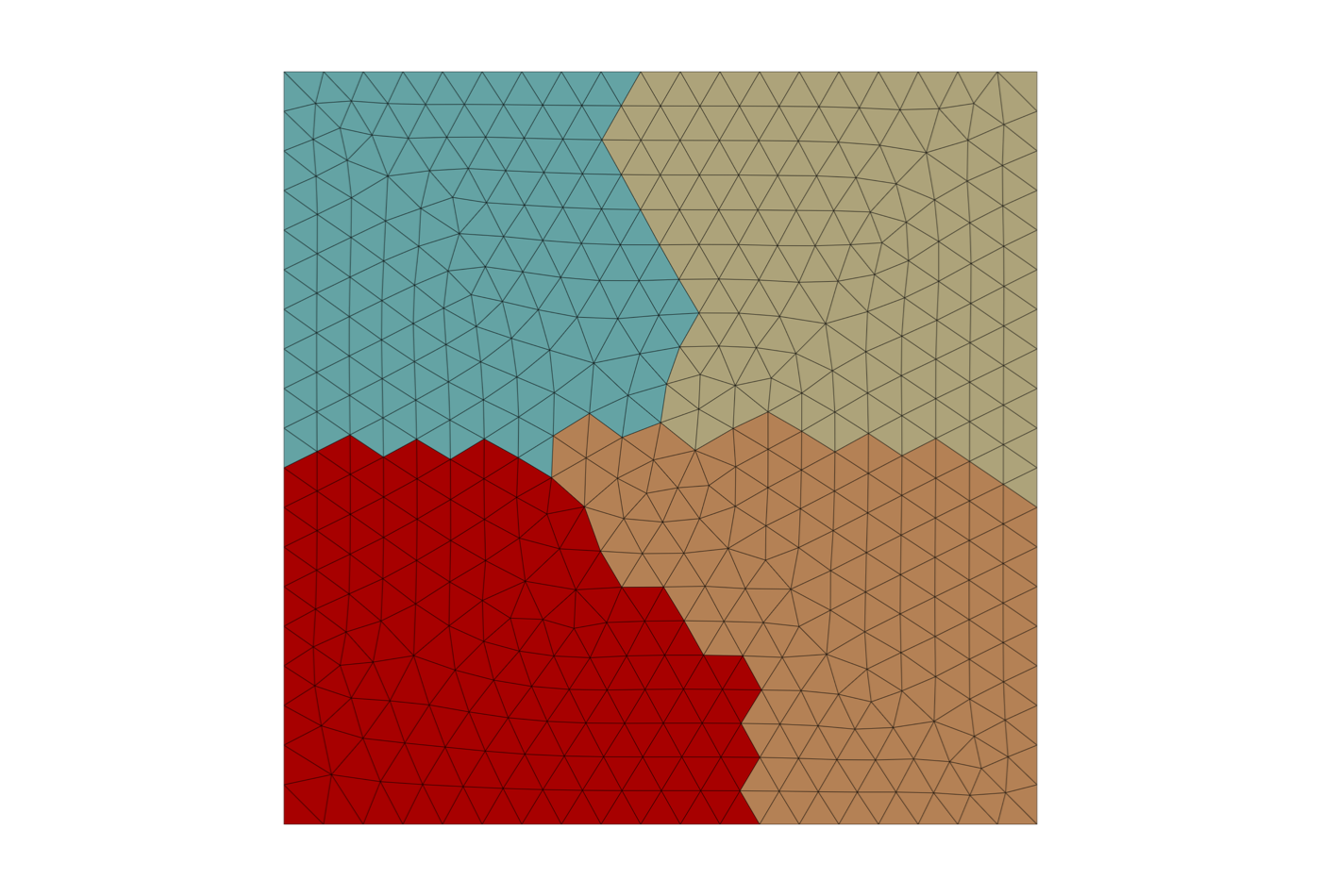}
\hfill
\maybeimage[width=.32\linewidth, trim={250 80 180 80}, clip]{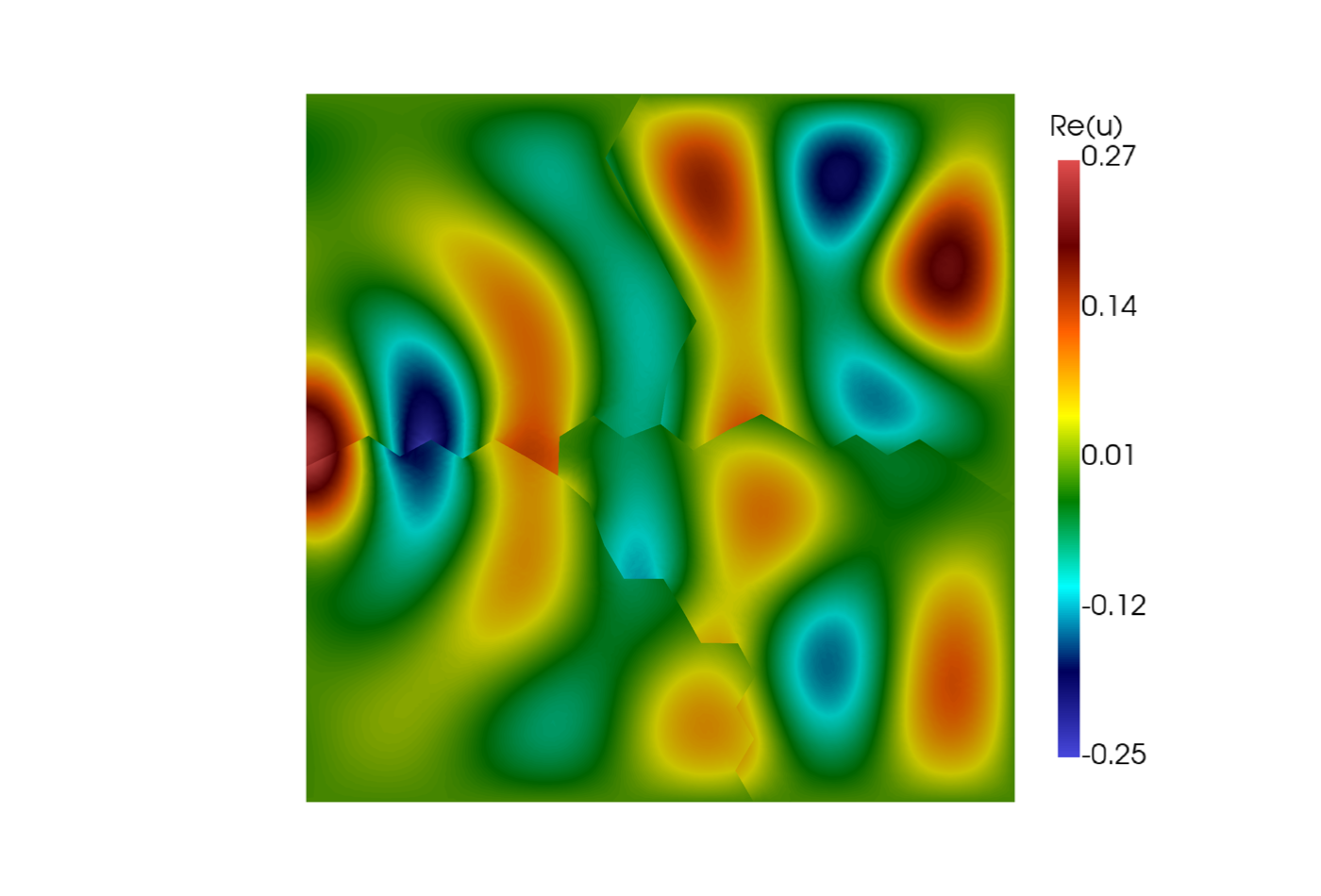}
\hfill
\maybeimage[width=.32\linewidth, trim={250 80 180 80}, clip]{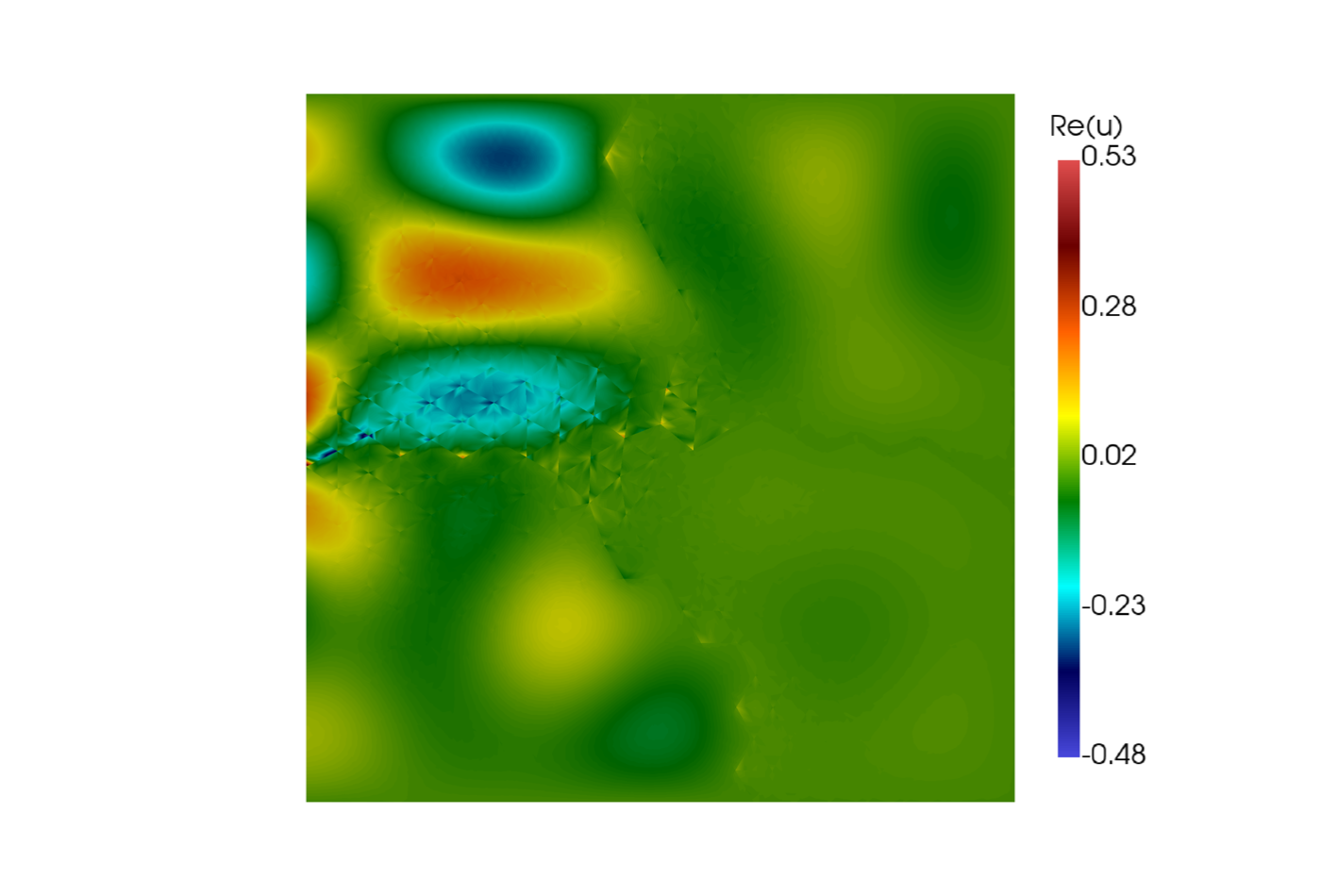}
\vspace{-0.5em}
\caption{
\textit{One application} of the multiplicative preconditioner on a unit square with a plane wave coming from the left.
Left: the four-subdomain METIS partition.
Middle: uses $\alpha = \beta = \delta = 1/2$ in the DG formulation.
Right: the same experiment with $\beta$ set to zero.
}
\label{fig:sweepsolver-visualization}
\end{figure}

To illustrate the importance of the normal-derivative jump term, 
we consider a simple example with an open propagation geometry.
In \Cref{fig:sweepsolver-visualization} we show a four-subdomain partition, 
the real part of the solution after one application of the multiplicative preconditioner with the standard DG coupling, 
and the same experiment with the interior $\beta$ term set to zero.
After one iteration, the intermediate iterate is still strongly shaped by the local subdomain solves.
Removing the normal-derivative jump term, i.e.\ setting $\beta=0$, produces stronger local reflections.

\section{Numerical experiments}
\label{sec:numerics}

Unless stated otherwise, timings were performed on a dual-socket machine with two Intel Xeon E5-2687W v4 processors (24 physical cores, 3.0 GHz), 512 GB of RAM, using 24 threads.
For the implementation of the methods we are using \texttt{NGSolve} \cite{ngsolve} and \texttt{NGSTrefftz} \cite{ngstrefftz}.
Replication data are available in \cite{GPSV26_data}.

\subsection{Parameter discussion}\label{sec:param-disc}
We first discuss the choice of penalty (stabilization) parameters. 
To that end, we consider the unit square with the exact plane wave solution
\[
    u_\star(x,y) := \exp(i\omega(\cos(\pi/5)x+\sin(\pi/5)y)),
    \qquad (x,y) \in \Omega,
\]
through the impedance datum
$g := \grad u_\star\cdot n+i\omega u_\star$ on the whole outer boundary. 
We fix the polynomial order $p=4$, a target resolution of five elements per wavelength, $\delta = 1/2$, 
and test six sets of parameters for both additive and multiplicative versions of the fixed point solver, 
Schwarz preconditioned CG (PCG), and Schwarz preconditioned GMRES (without restart) with \(J=16\) subdomains.
For all iterative methods we use the relative unpreconditioned residual in the stopping criteria, for this experiment the tolerance is set to \(10^{-6}\).
Since including the partition of unity (weighting) in the preconditioner breaks its symmetry, it is no longer applicable together with CG.

In \cite{Leumuller2023} the choice $\alpha = \beta = 1$ was proposed for HDG, 
while \cite{Leumuller2025} suggests $\alpha = 1/2, \beta = 2$. 
Furthermore, we consider $\alpha = \beta = 1/2$ as it recovers the ultra weak variational formulation \cite{CD98,BM08}.
We call parameters weighted, or weighted versions, if they are multiplied by \((p / \log(p)) / (h \omega)\) and \(h \omega / (p / \log(p))\), respectively, for \(\alpha\) and \(\beta\), as advocated for in~\cite{HMP11,M11}.

\Cref{tab:helmholtz-param-disc} reports iteration counts for all combinations of methods, penalty parameters, and preconditioner types.
Rows report methods, columns show the tuples $(\alpha,\beta)$, where the values indicated for \(\alpha\) and \(\beta\), unweighted or weighted, are multiplied by the factors above.

The lack of convergence for the fixed point solver with an additive preconditioner that uses overlap without a partition of unity is to be expected due to the over-correction in the overlap.
On the contrary, we observe that, in the multiplicative version, a partition of unity is not necessary since the approximation is updated sequentially.
Similarly, the sequential nature of multiplicative preconditioners implies that they typically require half the number of iterations of their additive counterparts (that are easily parallelized on the other hand).

The number of iterations for PCG and GMRES with the same preconditioner are comparable, so that PCG should be the preferred method to choose due to its short recurrence (reduced memory requirements compared to GMRES).
Overall we remark that Krylov-based methods are quite robust with respect to the choice of parameters.
With overlap, DG and TDG give the lowest multiplicative iteration counts.
Finally, we point out that results for a non-preconditioned CG are not reported because of a lack of convergence.

\addtolength{\tabcolsep}{-0.25em}
\newcommand{\ParamDiscTable}[2]{%
\begin{minipage}{0.49\linewidth}
\centering
\vspace{0.01\textheight}%
\small{#1}\\[0.3em]
\begingroup
\resizebox{\linewidth}{!}{%
\pgfplotstabletypeset[
    col sep=comma,
    string type,
    columns={
	method,overlap,pou,a0_b0,a0_b2,a1_b1,am0_bm0,am0_bm2,am1_bm1
    }, 
    columns/method/.style={column name={method}, column type=c},
    columns/overlap/.style={column name={overlap}, column type=c},
    columns/pou/.style={column name={PoU}, column type=c},
    columns/a0_b0/.style={column name={$(\nicefrac{1}{2},\nicefrac{1}{2})$}, column type=c},
    columns/a0_b2/.style={column name={$(\nicefrac{1}{2},{2})$}, column type=c},
    columns/a1_b1/.style={column name={$(1,1)$}, column type=c},
    columns/am0_bm0/.style={column name={$(\nicefrac{1}{2},\nicefrac{1}{2})$}, column type=c},
    columns/am0_bm2/.style={column name={$(\nicefrac{1}{2},{2})$}, column type=c},
    columns/am1_bm1/.style={column name={$(1,1)$}, column type=c},
    every head row/.style={
        before row=\toprule,
        after row=\midrule,
    },
    every head row/.style={
        before row={
            \toprule
            & & & \multicolumn{3}{c}{constant parameters}
            & \multicolumn{3}{c}{weighted parameters}\\
            \cmidrule(lr){4-6}
            \cmidrule(lr){7-9}
        },
    after row=\midrule,
},
    every odd row/.style={
        before row={\rowcolor{gray!15}}
    },
    every last row/.style={after row=\bottomrule},
]{#2}
\vspace{0.01\textheight}%
}%
\endgroup
\end{minipage}}

\begin{table}[ht!]
\centering
\ParamDiscTable{Additive fixed point}{./numex/unit_square_table_param_disc_jacobi_fixed.csv}
\ParamDiscTable{Multiplicative fixed point}{./numex/unit_square_table_param_disc_gs_fixed.csv}

\ParamDiscTable{Additive PCG}{./numex/unit_square_table_param_disc_jacobi_cg.csv}
\ParamDiscTable{Multiplicative PCG}{./numex/unit_square_table_param_disc_gs_cg.csv}

\ParamDiscTable{Additive GMRES}{./numex/unit_square_table_param_disc_jacobi_gmres.csv}
\ParamDiscTable{Multiplicative GMRES}{./numex/unit_square_table_param_disc_gs_gmres.csv}

\ParamDiscTable{Unpreconditioned GMRES}{./numex/unit_square_table_param_disc_gmres.csv}

\vspace{0.5em}
\caption{
Iteration counts on the unit square with exact plane wave solution for different methods, solvers and parameters \((\alpha,\beta)\). 
Rows report methods, columns show the tuples $(\alpha,\beta)$, where the values indicated for \(\alpha\) and \(\beta\) 
are multiplied by \((p / \log(p)) / (h \omega)\) and \(h \omega / (p / \log(p))\) respectively in the weighted case.
}
\label{tab:helmholtz-param-disc}
\vspace{-1em}
\end{table}

Concluding, together with the results for the non-preconditioned GMRES solver, this highlights that preconditioners are mandatory to obtain efficient solvers in this context (simple fixed point solvers with additive preconditioners require less iterations).
The parameter shows that Krylov acceleration makes the methods comparatively insensitive to the tested parameter variations.
Based on all these results, and further testing in other experiments, we from now on fix the constant parameters $\alpha = \beta = 1$ for HDG, and $\alpha=\beta=\delta=1/2$ for DG and TDG.

\subsection{Influence of \texorpdfstring{$h$}{h} and \texorpdfstring{$p$}{p} refinement}

We consider the exact solution on the unit square as in \Cref{sec:param-disc}.
And an analogous unit-cube test case in three dimensions with plane-wave direction $(1,2,3)/\sqrt{14}$.
We are interested in the influence of $h$- and $p$-refinement on the convergence of the solver.
We compare the discretizations using CG with a multiplicative preconditioner.

In the $h$-refinement case, the polynomial order $p = 5$ is kept fixed, and the number of elements per wavelength ranges from roughly $2^{i/2}$, $i=1,\dots,6$.
In the $p$-refinement case, the target resolution is kept fixed at five elements per wavelength, and the polynomial order ranges from $p=2$ to $p=7$.
The convergence tolerance is \(10^{-12}\) for this experiment.

\begin{figure}[ht!]
\centering
\resizebox{\linewidth}{!}{
\begin{tikzpicture}
    \begin{groupplot}[
        group style={
            group size=3 by 2,
            horizontal sep=1.4cm,
            vertical sep=1.4cm,
        },
        width=.52\linewidth,
        height=.45\linewidth,
        ymajorgrids=true,
        xmajorgrids=true,
        grid style=dashed,
        unbounded coords=discard,
        filter discard warning=false,
        cycle list name=paulcolors6,
        x dir=reverse,
        legend style={
            at={(0.97,0.97)},
            anchor=north east,
            draw=none,
            fill=none,
            font=\normalsize,
            legend style={at={(1.8, -1.65)}, anchor=north, draw=none, fill=none, legend columns=6},
        },
    ]
    \nextgroupplot[
        xlabel={mesh size $h$},
        ylabel={$L^2$ error},
        xmin=0.068,
        xmax=0.7,
        xmode=log,
        ymode=log,
    ]
    \addplot+[discard if not={overlap}{0}]
        table [x=mesh_size, y=l2_error, col sep=comma]
        {./numex/helmholtz_unit_square_dg_h.csv};
    \addplot+[discard if not={overlap}{0}]
        table [x=mesh_size, y=l2_error, col sep=comma]
        {./numex/helmholtz_unit_square_tdg_h.csv};
    \addplot+[]
        table [x=mesh_size, y=l2_error, col sep=comma]
        {./numex/helmholtz_unit_square_hdg_h.csv};
    \addplot+[discard if not={overlap}{1}]
        table [x=mesh_size, y=l2_error, col sep=comma]
        {./numex/helmholtz_unit_square_dg_h.csv};
    \addplot+[discard if not={overlap}{1}]
        table [x=mesh_size, y=l2_error, col sep=comma]
        {./numex/helmholtz_unit_square_tdg_h.csv};
        \legend{DG (without overlap),TDG (without overlap),HDG,DG,TDG}
    \nextgroupplot[
        xlabel={mesh size $h$},
        ylabel={iterations},
        xmin=0.068,
        xmax=0.7,
        xmode=log,
    ]
    \addplot+[discard if not={overlap}{0}]
        table [x=mesh_size, y=steps, col sep=comma]
        {./numex/helmholtz_unit_square_dg_h.csv};
    \addplot+[discard if not={overlap}{0}]
        table [x=mesh_size, y=steps, col sep=comma]
        {./numex/helmholtz_unit_square_tdg_h.csv};
    \addplot+[]
        table [x=mesh_size, y=steps, col sep=comma]
        {./numex/helmholtz_unit_square_hdg_h.csv};
    \addplot+[discard if not={overlap}{1}]
        table [x=mesh_size, y=steps, col sep=comma]
        {./numex/helmholtz_unit_square_dg_h.csv};
    \addplot+[discard if not={overlap}{1}]
        table [x=mesh_size, y=steps, col sep=comma]
        {./numex/helmholtz_unit_square_tdg_h.csv};
    \nextgroupplot[
        xlabel={mesh size $h$},
        ylabel={avg. local nnze},
        xmode=log,
        ymode=log,
        xmin=0.068,
        xmax=0.7,
    ]
    \addplot+[discard if not={overlap}{0}]
        table [x=mesh_size, y=avg_domain_nnze, col sep=comma]
        {./numex/helmholtz_unit_square_dg_h.csv};
    \addplot+[discard if not={overlap}{0}]
        table [x=mesh_size, y=avg_domain_nnze, col sep=comma]
        {./numex/helmholtz_unit_square_tdg_h.csv};
    \addplot+[]
        table [x=mesh_size, y=avg_domain_nnze, col sep=comma]
        {./numex/helmholtz_unit_square_hdg_h.csv};
    \addplot+[discard if not={overlap}{1}]
        table [x=mesh_size, y=avg_domain_nnze, col sep=comma]
        {./numex/helmholtz_unit_square_dg_h.csv};
    \addplot+[discard if not={overlap}{1}]
        table [x=mesh_size, y=avg_domain_nnze, col sep=comma]
        {./numex/helmholtz_unit_square_tdg_h.csv};
    \nextgroupplot[
        xlabel={mesh size $h$},
        ylabel={$L^2$ error},
        xmin=0.068,
        xmax=0.7,
        xmode=log,
        ymode=log,
    ]
    \addplot+[discard if not={overlap}{0}]
        table [x=mesh_size, y=l2_error, col sep=comma]
        {./numex/helmholtz_unit_cube_dg_h.csv};
    \addplot+[discard if not={overlap}{0}]
        table [x=mesh_size, y=l2_error, col sep=comma]
        {./numex/helmholtz_unit_cube_tdg_h.csv};
    \addplot+[]
        table [x=mesh_size, y=l2_error, col sep=comma]
        {./numex/helmholtz_unit_cube_hdg_h.csv};
    \addplot+[discard if not={overlap}{1}]
        table [x=mesh_size, y=l2_error, col sep=comma]
        {./numex/helmholtz_unit_cube_dg_h.csv};
    \addplot+[discard if not={overlap}{1}]
        table [x=mesh_size, y=l2_error, col sep=comma]
        {./numex/helmholtz_unit_cube_tdg_h.csv};
    \nextgroupplot[
        xlabel={mesh size $h$},
        ylabel={iterations},
        xmin=0.068,
        xmax=0.7,
        xmode=log,
    ]
    \addplot+[discard if not={overlap}{0}]
        table [x=mesh_size, y=steps, col sep=comma]
        {./numex/helmholtz_unit_cube_dg_h.csv};
    \addplot+[discard if not={overlap}{0}]
        table [x=mesh_size, y=steps, col sep=comma]
        {./numex/helmholtz_unit_cube_tdg_h.csv};
    \addplot+[]
        table [x=mesh_size, y=steps, col sep=comma]
        {./numex/helmholtz_unit_cube_hdg_h.csv};
    \addplot+[discard if not={overlap}{1}]
        table [x=mesh_size, y=steps, col sep=comma]
        {./numex/helmholtz_unit_cube_dg_h.csv};
    \addplot+[discard if not={overlap}{1}]
        table [x=mesh_size, y=steps, col sep=comma]
        {./numex/helmholtz_unit_cube_tdg_h.csv};
    \nextgroupplot[
        xlabel={mesh size $h$},
        ylabel={avg. local nnze},
        xmode=log,
        ymode=log,
        xmin=0.068,
        xmax=0.7,
    ]
    \addplot+[discard if not={overlap}{0}]
        table [x=mesh_size, y=avg_domain_nnze, col sep=comma]
        {./numex/helmholtz_unit_cube_dg_h.csv};
    \addplot+[discard if not={overlap}{0}]
        table [x=mesh_size, y=avg_domain_nnze, col sep=comma]
        {./numex/helmholtz_unit_cube_tdg_h.csv};
    \addplot+[]
        table [x=mesh_size, y=avg_domain_nnze, col sep=comma]
        {./numex/helmholtz_unit_cube_hdg_h.csv};
    \addplot+[discard if not={overlap}{1}]
        table [x=mesh_size, y=avg_domain_nnze, col sep=comma]
        {./numex/helmholtz_unit_cube_dg_h.csv};
    \addplot+[discard if not={overlap}{1}]
        table [x=mesh_size, y=avg_domain_nnze, col sep=comma]
        {./numex/helmholtz_unit_cube_tdg_h.csv};
    \end{groupplot}
\end{tikzpicture}}
\vspace{-1em}
\caption{
Influence of $h$ refinement on DG, TDG and HDG on an exact plane wave
in 2D (top) and 3D (bottom).
When available results with non-overlapping partitions are included.
Left: $L^2$ error versus mesh size $h$.
Middle: iterations versus mesh size $h$.
Right: average non-zero entries in local matrix $A_{j}$ versus mesh size $h$.
}
\label{fig:helmholtz-exact-convergence-h}
\end{figure}
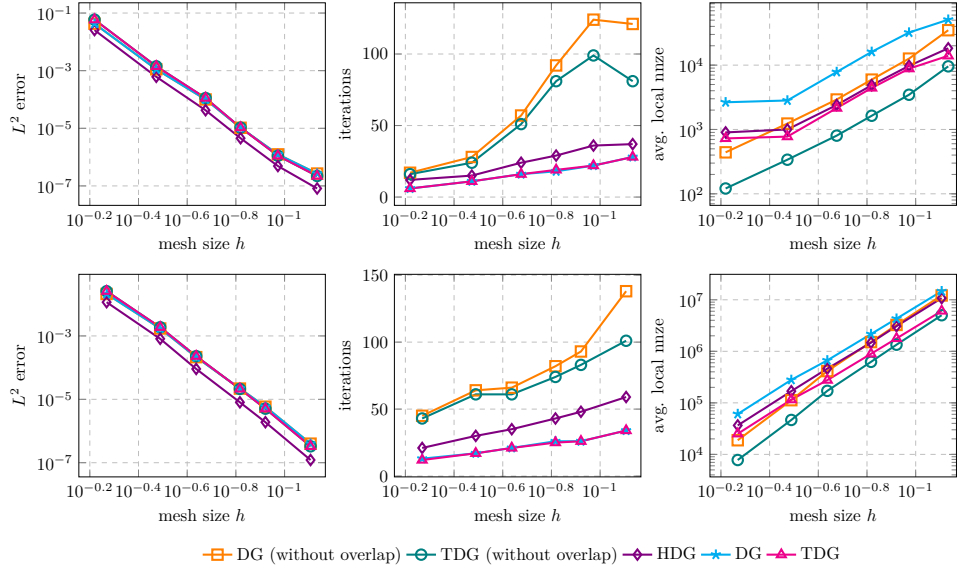
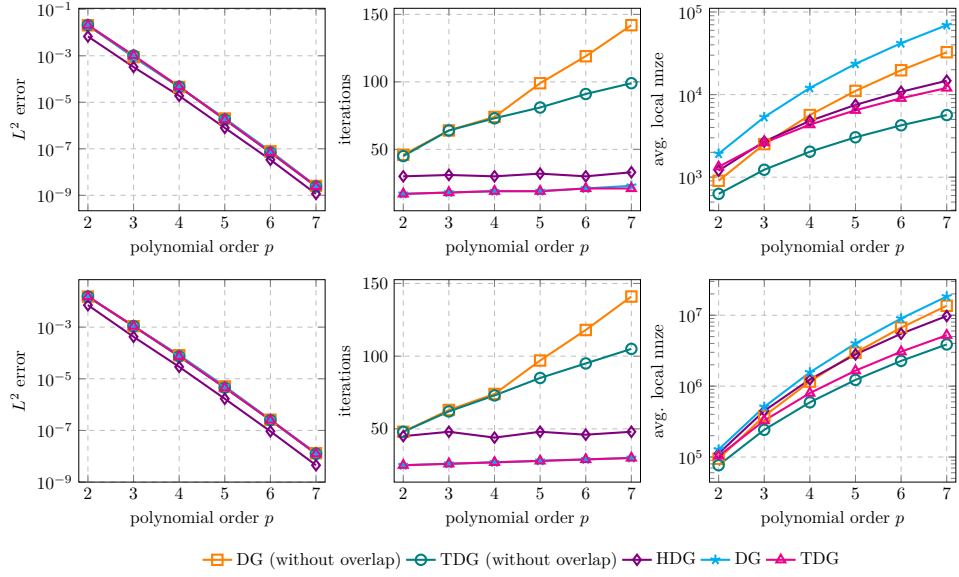
\begin{figure}[ht!]
\centering
\resizebox{\linewidth}{!}{
\begin{tikzpicture}
    \begin{groupplot}[
        group style={
            group size=3 by 2,
            horizontal sep=1.4cm,
            vertical sep=1.4cm,
        },
        width=.52\linewidth,
        height=.45\linewidth,
        ymajorgrids=true,
        xmajorgrids=true,
        grid style=dashed,
        unbounded coords=discard,
        filter discard warning=false,
        cycle list name=paulcolors6,
        legend style={
            at={(0.97,0.97)},
            anchor=north east,
            draw=none,
            fill=none,
            font=\normalsize,
            legend style={at={(1.8, -1.65)}, anchor=north, draw=none, fill=none, legend columns=6},
        },
    ]
    \nextgroupplot[
        xlabel={polynomial order $p$},
        ylabel={$L^2$ error},
        xmin=1.8,
        xmax=7.2,
        ymode=log,
    ]
    \addplot+[discard if not={overlap}{0}]
        table [x=order, y=l2_error, col sep=comma]
        {./numex/helmholtz_unit_square_dg.csv};
    \addplot+[discard if not={overlap}{0}]
        table [x=order, y=l2_error, col sep=comma]
        {./numex/helmholtz_unit_square_tdg.csv};
    \addplot+[]
        table [x=order, y=l2_error, col sep=comma]
        {./numex/helmholtz_unit_square_hdg.csv};
    \addplot+[discard if not={overlap}{1}]
        table [x=order, y=l2_error, col sep=comma]
        {./numex/helmholtz_unit_square_dg.csv};
    \addplot+[discard if not={overlap}{1}]
        table [x=order, y=l2_error, col sep=comma]
        {./numex/helmholtz_unit_square_tdg.csv};
        \legend{DG (without overlap),TDG (without overlap),HDG,DG,TDG}
    \nextgroupplot[
        xlabel={polynomial order $p$},
        ylabel={iterations},
        xmin=1.8,
        xmax=7.2,
        xtick=data,
    ]
    \addplot+[discard if not={overlap}{0}]
        table [x=order, y=steps, col sep=comma]
        {./numex/helmholtz_unit_square_dg.csv};
    \addplot+[discard if not={overlap}{0}]
        table [x=order, y=steps, col sep=comma]
        {./numex/helmholtz_unit_square_tdg.csv};
    \addplot+[]
        table [x=order, y=steps, col sep=comma]
        {./numex/helmholtz_unit_square_hdg.csv};
    \addplot+[discard if not={overlap}{1}]
        table [x=order, y=steps, col sep=comma]
        {./numex/helmholtz_unit_square_dg.csv};
    \addplot+[discard if not={overlap}{1}]
        table [x=order, y=steps, col sep=comma]
        {./numex/helmholtz_unit_square_tdg.csv};
    \nextgroupplot[
        xlabel={polynomial order $p$},
        ylabel={avg. local nnze},
        ymode=log,
        xmin=1.8,
        xmax=7.2,
        xtick=data,
    ]
    \addplot+[discard if not={overlap}{0}]
        table [x=order, y=avg_domain_nnze, col sep=comma]
        {./numex/helmholtz_unit_square_dg.csv};
    \addplot+[discard if not={overlap}{0}]
        table [x=order, y=avg_domain_nnze, col sep=comma]
        {./numex/helmholtz_unit_square_tdg.csv};
    \addplot+[]
        table [x=order, y=avg_domain_nnze, col sep=comma]
        {./numex/helmholtz_unit_square_hdg.csv};
    \addplot+[discard if not={overlap}{1}]
        table [x=order, y=avg_domain_nnze, col sep=comma]
        {./numex/helmholtz_unit_square_dg.csv};
    \addplot+[discard if not={overlap}{1}]
        table [x=order, y=avg_domain_nnze, col sep=comma]
        {./numex/helmholtz_unit_square_tdg.csv};
    \nextgroupplot[
        xlabel={polynomial order $p$},
        ylabel={$L^2$ error},
        xmin=1.8,
        xmax=7.2,
        ymode=log,
    ]
    \addplot+[]
        table [x=order, y=l2_error, col sep=comma]
        {./numex/helmholtz_unit_cube_dg.csv};
    \addplot+[]
        table [x=order, y=l2_error, col sep=comma]
        {./numex/helmholtz_unit_cube_tdg.csv};
    \addplot+[]
        table [x=order, y=l2_error, col sep=comma]
        {./numex/helmholtz_unit_cube_hdg.csv};
    \addplot+[]
        table [x=order, y=l2_error, col sep=comma]
        {./numex/helmholtz_unit_cube_dg_overlap.csv};
    \addplot+[]
        table [x=order, y=l2_error, col sep=comma]
        {./numex/helmholtz_unit_cube_tdg_overlap.csv};
    \nextgroupplot[
        xlabel={polynomial order $p$},
        ylabel={iterations},
        xmin=1.8,
        xmax=7.2,
        xtick=data,
    ]
    \addplot+[]
        table [x=order, y=steps, col sep=comma]
        {./numex/helmholtz_unit_cube_dg.csv};
    \addplot+[]
        table [x=order, y=steps, col sep=comma]
        {./numex/helmholtz_unit_cube_tdg.csv};
    \addplot+[]
        table [x=order, y=steps, col sep=comma]
        {./numex/helmholtz_unit_cube_hdg.csv};
    \addplot+[]
        table [x=order, y=steps, col sep=comma]
        {./numex/helmholtz_unit_cube_dg_overlap.csv};
    \addplot+[]
        table [x=order, y=steps, col sep=comma]
        {./numex/helmholtz_unit_cube_tdg_overlap.csv};
    \nextgroupplot[
        xlabel={polynomial order $p$},
        ylabel={avg. local nnze},
        ymode=log,
        xmin=1.8,
        xmax=7.2,
        xtick=data,
    ]
    \addplot+[]
        table [x=order, y=avg_domain_nnze, col sep=comma]
        {./numex/helmholtz_unit_cube_dg.csv};
    \addplot+[]
        table [x=order, y=avg_domain_nnze, col sep=comma]
        {./numex/helmholtz_unit_cube_tdg.csv};
    \addplot+[]
        table [x=order, y=avg_domain_nnze, col sep=comma]
        {./numex/helmholtz_unit_cube_hdg.csv};
    \addplot+[]
        table [x=order, y=avg_domain_nnze, col sep=comma]
        {./numex/helmholtz_unit_cube_dg_overlap.csv};
    \addplot+[]
        table [x=order, y=avg_domain_nnze, col sep=comma]
        {./numex/helmholtz_unit_cube_tdg_overlap.csv};
    \end{groupplot}
\end{tikzpicture}}
\vspace{-1em}
\caption{
Influence of $p$ refinement on DG, TDG and HDG on an exact plane wave
in 2D (top) and 3D (bottom).
When available results with non-overlapping partitions are included.
Left: $L^2$ error versus polynomial order.
Middle: iterations versus polynomial order.
Right: average non-zero entries in local matrix $A_{j}$ versus polynomial order.
}
\label{fig:helmholtz-exact-convergence}
\end{figure}

\Cref{fig:helmholtz-exact-convergence-h} and~\Cref{fig:helmholtz-exact-convergence} depict $L^2$-error, iteration count and number of non-zeros in local matrices $A_{j}$ with respect to mesh size $h$ and polynomial order \(p\) respectively. 
The $L^2$-errors are very similar for all methods (with a slight advantage for HDG), which makes the comparison in iteration count and matrix size valid.
Iteration counts for non-overlapping DG and TDG are not robust with respect to polynomial order, while HDG, overlapping DG and TDG show more robust behaviour, particularly in the case of $p$-refinement.
The robustness with respect to order in overlapping DD methods with impedance boundary condition is a known phenomenon, see~\cite{Gong2023}.

The smallest iterations counts appear to be obtained for DG and TDG with overlap, although iteration counts for HDG are competitive.
As expected, the local matrix sizes increase with the polynomial degree, the ones for DG having the largest relative and absolute increase, compare \Cref{tab:helmholtz-simplicial-counts}.
We reported this quantity as an indicator of the computational effort required to solve the local problems.
HDG and TDG have a smaller increase in local matrix size as they reduce the number of degrees of freedom from $\mathcal{O}(p^d)$ to $\mathcal{O}(p^{d-1})$ in $d$ dimensions.

\begin{table}[ht!]
\centering
\footnotesize
\setlength{\tabcolsep}{0.25em}
\begin{filecontents*}{./numex/helmholtz_table_timings_p7.csv}
method,overlap,pou,assembly_time,setup_time,solve_time,avgmatvec_time,postprocess_time,total_time,steps,avg_domain_nnze
DG,\xmark,\xmark,17.862906977534294,252.95761902444065,126.57260755822062,0.9171928083929031,3.0891969799995422e-06,397.6772736189887,139,13616100.0
TDG,\xmark,\xmark,20.20257391408086,68.42953842971474,35.33983619790524,0.3398061172875504,0.019389980472624302,124.26893303263932,105,3873024.0
HDG,\cmark,\xmark,67.01760153099895,168.50780819077045,33.29912842810154,0.723894096263077,2.8765770960599184,271.9994188202545,47,9754332.0
DG,\cmark,\xmark,17.66625460051,280.52891249395907,38.24070646800101,1.3186450506207243,1.6726553440093994e-06,336.72192649077624,30,18368100.0
TDG,\cmark,\xmark,20.03222766518593,79.99571749009192,12.82177973818034,0.44213033579932204,0.02533307857811451,113.15392081439495,30,5224704.0
\end{filecontents*}
\pgfplotstabletypeset[
    col sep=comma, fixed, precision=1, fixed zerofill=true,
    columns={
    method,overlap,pou,assembly_time,setup_time,avgmatvec_time,solve_time,postprocess_time,total_time,steps
}, 
    columns/method/.style={column name={method}, column type=c, string type},
    columns/overlap/.style={column name={overlap}, column type=c, string type},
    columns/pou/.style={column name={PoU}, column type=c, string type},
    columns/assembly_time/.style={column name={assembly}, column type=c},
    columns/setup_time/.style={column name={factorization}, column type=c},
    columns/avgmatvec_time/.style={column name={avg. it.}, column type=c},
    columns/solve_time/.style={column name={solve}, column type=c},
    columns/postprocess_time/.style={column name={postprocess}, column type=c},
    columns/total_time/.style={column name={total}, column type=c},
    columns/steps/.style={column name={iterations}, column type=c, string type},
    every head row/.style={
        before row=\toprule,
        after row=\midrule,
    },
    every last row/.style={after row=\bottomrule},
    every odd row/.style={
        before row={\rowcolor{gray!15}}
    },
    every last row/.style={after row=\bottomrule},
]{./numex/helmholtz_table_timings_p7.csv}
\vspace{.5em}
\caption{
Timings for CG with multiplicative preconditioner in seconds on the unit cube with exact plane wave solution for polynomial order $p=7$. 
}
\label{tab:helmholtz-timings}
\vspace{-1em}
\end{table}

As runtimes become more significant in three dimensions, \Cref{tab:helmholtz-timings} additionally 
reports wall time for $p=7$ and five elements per wavelength.
We compare timings for matrix assembly (including assembly of the embedded Trefftz basis functions and elimination in HDG), factorization of local matrices for the preconditioner, average iteration cost, total solving, postprocessing (HDG only), and iteration count. 
It appears the computation of the embedded Trefftz basis functions in TDG does not increase significantly the assembly time compared to DG (for which larger matrices have to be assembled).
The assembly of HDG is significantly larger despite the smaller condensed matrix due to the elimination step.
The factorization time is a function of the size of the local matrices.
The solving time reflects iteration count and average iteration cost.

\subsection{Strong and weak scaling tests}

We continue with scaling tests on the same problem on the unit square, now with $\omega=40$, where we study the influence of the number of subdomains \(J\) on the performance of the preconditioners, for the multiplicative variant only, using the CG solver.
We provide results for a \textit{strong scaling test}, namely we fix all parameters and increase the number of subdomains \(J\) from \(2\) to \(64\).
We also provide the results for a \textit{weak scaling test} in which we fix the geometry and increase both \(\omega\) and \(J\) such that the number of degrees of freedom per subdomain stays approximately constant: $\omega$ ranges from $5$ to $80$ and \(J\) from $5$ to $1280$.

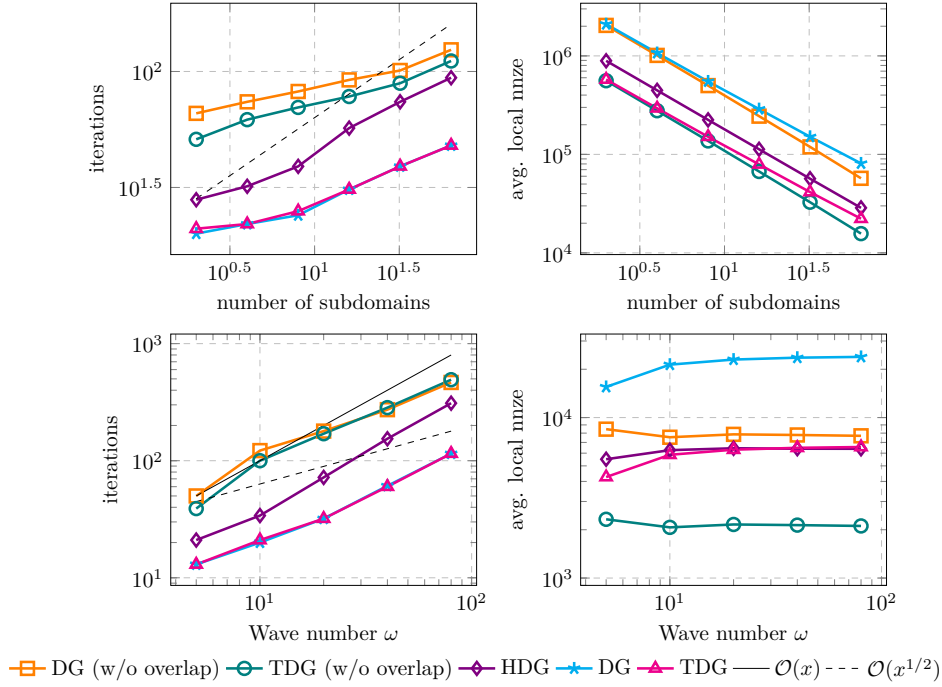
\begin{figure}[ht!]
\centering
\resizebox{\linewidth}{!}{
\begin{tikzpicture}
    \begin{groupplot}[
        group style={
            group size=2 by 2,
            horizontal sep=1.7cm,
            vertical sep=1.3cm,
        },
        width=.52\linewidth,
        height=.45\linewidth,
        ymajorgrids=true,
        xmajorgrids=true,
        grid style=dashed,
        unbounded coords=discard,
        filter discard warning=false,
        cycle list name=paulcolors6,
        legend style={
            at={(0.97,0.97)},
            anchor=north east,
            draw=none,
            fill=none,
            font=\normalsize,
            legend style={at={(1.0, -0.25)}, anchor=north, draw=none, fill=none, legend columns=7},
        },
    ]
    \nextgroupplot[
        xlabel={number of subdomains},
        ylabel={iterations},
        xmode=log,
        ymode=log,
    ]
    \addplot+[discard if not={overlap}{0}]
        table [x=ndomains, y=steps, col sep=comma]
        {./numex/helmholtz_unit_square_dg_subdomains.csv};
    \addplot+[discard if not={overlap}{0}]
        table [x=ndomains, y=steps, col sep=comma]
        {./numex/helmholtz_unit_square_tdg_subdomains.csv};
    \addplot+[]
        table [x=ndomains, y=steps, col sep=comma]
        {./numex/helmholtz_unit_square_hdg_subdomains.csv};
    \addplot+[discard if not={overlap}{1}]
        table [x=ndomains, y=steps, col sep=comma]
        {./numex/helmholtz_unit_square_dg_subdomains.csv};
    \addplot+[discard if not={overlap}{1}]
        table [x=ndomains, y=steps, col sep=comma]
        {./numex/helmholtz_unit_square_tdg_subdomains.csv};
          \addplot[domain=2:64,dashed] {20*x^(0.5)};
    \nextgroupplot[
        xlabel={number of subdomains},
        ylabel={avg. local nnze},
        xmode=log,
        ymode=log,
    ]
    \addplot+[discard if not={overlap}{0}]
        table [x=ndomains, y=avg_domain_nnze, col sep=comma]
        {./numex/helmholtz_unit_square_dg_subdomains.csv};
    \addplot+[discard if not={overlap}{0}]
        table [x=ndomains, y=avg_domain_nnze, col sep=comma]
        {./numex/helmholtz_unit_square_tdg_subdomains.csv};
    \addplot+[]
        table [x=ndomains, y=avg_domain_nnze, col sep=comma]
        {./numex/helmholtz_unit_square_hdg_subdomains.csv};
    \addplot+[discard if not={overlap}{1}]
        table [x=ndomains, y=avg_domain_nnze, col sep=comma]
        {./numex/helmholtz_unit_square_dg_subdomains.csv};
    \addplot+[discard if not={overlap}{1}]
        table [x=ndomains, y=avg_domain_nnze, col sep=comma]
        {./numex/helmholtz_unit_square_tdg_subdomains.csv};
        
    \nextgroupplot[
        xlabel={Wave number $\omega$},
        ylabel={iterations},
        xmode=log,
        ymode=log,
    ]
    \addplot+[discard if not={overlap}{0}]
        table [x=omega, y=steps, col sep=comma]
        {./numex/helmholtz_unit_square_dg_iomega.csv};
    \addplot+[discard if not={overlap}{0}]
        table [x=omega, y=steps, col sep=comma]
        {./numex/helmholtz_unit_square_tdg_iomega.csv};
    \addplot+[]
        table [x=omega, y=steps, col sep=comma]
        {./numex/helmholtz_unit_square_hdg_iomega.csv};
    \addplot+[discard if not={overlap}{1}]
        table [x=omega, y=steps, col sep=comma]
        {./numex/helmholtz_unit_square_dg_iomega.csv};
    \addplot+[discard if not={overlap}{1}]
        table [x=omega, y=steps, col sep=comma]
        {./numex/helmholtz_unit_square_tdg_iomega.csv};
          \addplot[domain=5:80] {10*x};
          \addplot[domain=5:80,dashed] {20*x^(0.5)};
          \addlegendimage{solid}
          \addlegendimage{dashed}
        \legend{DG (w/o overlap),TDG (w/o overlap),HDG,DG,TDG, $\mathcal O(x)$, $\mathcal O(x^{1/2})$}
    \nextgroupplot[
        xlabel={Wave number $\omega$},
        ylabel={avg. local nnze},
        xmode=log,
        ymode=log,
        ymin=900
    ]
    \addplot+[discard if not={overlap}{0}]
        table [x=omega, y=avg_domain_nnze, col sep=comma]
        {./numex/helmholtz_unit_square_dg_iomega.csv};
    \addplot+[discard if not={overlap}{0}]
        table [x=omega, y=avg_domain_nnze, col sep=comma]
        {./numex/helmholtz_unit_square_tdg_iomega.csv};
    \addplot+[]
        table [x=omega, y=avg_domain_nnze, col sep=comma]
        {./numex/helmholtz_unit_square_hdg_iomega.csv};
    \addplot+[discard if not={overlap}{1}]
        table [x=omega, y=avg_domain_nnze, col sep=comma]
        {./numex/helmholtz_unit_square_dg_iomega.csv};
    \addplot+[discard if not={overlap}{1}]
        table [x=omega, y=avg_domain_nnze, col sep=comma]
        {./numex/helmholtz_unit_square_tdg_iomega.csv};
    \end{groupplot}
\end{tikzpicture}}
\vspace{-1em}
\caption{
Strong (top) and weak (bottom) scaling tests for DG, TDG and HDG on an exact plane wave
in 2D.
When available results with non-overlapping partitions are included.
Left: iteration versus number of subdomains.
Right: average non-zero entries in local matrix $A_{j}$ versus number of subdomains.
}
\label{fig:helmholtz-exact-subdomains}
\end{figure}

The results reported in~\Cref{fig:helmholtz-exact-subdomains} show the expected behavior: One-level domain decomposition preconditioners are not robust with respect to the number of subdomains.
For additive type preconditioners, which only transfer information between neighboring subdomains at each iteration, we expect a square-root increase in the iteration count in the strong scaling test.
For the (forward and backward) multiplicative preconditioners considered here, global transfer of information remains partial, and a sub-linear increase in the number of iterations is observed in both scaling tests.

The average local number of non-zeros in the local matrices \(A_{j}\) is also reported, and agrees with the partitioning strategy: in the strong scaling test, the average number of non-zeros decreases with increasing number of subdomains, while in the weak scaling test it remains approximately constant.

\subsection{Trapping geometries and solver comparison}

We consider three two-dimensional trapping geometries in the unit square.
They are bounded-domain analogues of the standard exterior examples for elliptic, parabolic, and hyperbolic trapping.
\Cref{fig:trapping-geometries} shows the three geometries and the relevant geometric parameters.
We consider the scattering of a plane wave (coming from the left) by the gray sets considered as Dirichlet obstacles.
The experiments below use $\omega=40$, six elements per wavelength and polynomial order $p=4$;
and the multiplicative preconditioners are used with the CG solver.

\input{trapping-geometries}

For the hyperbolic case, two disjoint smooth strictly convex obstacles form the
standard Ikawa-type example; the trapped ray between the obstacles is geometric
and does not require a special wavenumber \cite{Ikawa1988,Burq2004,
ChandlerWildeSpenceGibbsSmyshlyaev2020}.
The precise geometric parameters are $x_h=0.62$, $D=0.44$ and $r=0.16$.

For the parabolic case, the relevant quantity is the interval condition
$\omega\in\{m\pi/a:m\in\IN\}$ for the distance $a$ between parallel Dirichlet
sides \cite{ChandlerWildeSpenceGibbsSmyshlyaev2020,
BetckeChandlerWildeGrahamLangdonLindner2011}.
The precise geometric parameters are $x_p=0.60$, $a=2\pi/\omega$, and $b=0.25$, giving the mode $m=2$ in $\omega=m\pi/a$.

For the elliptic case, we use the small elliptic cavity of
\cite{BetckeChandlerWildeGrahamLangdonLindner2011,MarchandGalkowskiSpenceSpence2021}:
the trapping frequencies are the Mathieu eigenfrequencies of the corresponding
ellipse, and the associated eigenfunctions localize near the minor-axis
bouncing-ball orbit.
The first eigenfrequency is given by $k_{1,0}\approx9.977$. 
The unscaled cavity has inner semi-axes $(1,1/2)$, outer semi-axes $(1.3,0.6)$, and $\phi=7\pi/10$.
We scale all the obstacle lengths by $s=k_{1,0}/40\approx0.2494$, so that the scaled cavity is tuned to $\omega=40$.
The precise geometric parameters are then $x_e=0.5$, $(a_i,b_i)\approx(0.2494,0.1247)$, $(a_o,b_o)\approx(0.3243,0.1497)$.

\begin{figure}[ht!]
\centering
\resizebox{\linewidth}{!}{
\begin{tikzpicture}
    \begin{groupplot}[
        group style={group size=3 by 3, horizontal sep=1.1cm, vertical sep=0.75cm},
        width=.47\linewidth,
        height=.31\linewidth,
        ymode=log,
        xmin=0,
        xmax=99,
        ymin=1e-9,
        ymax=1e3,
        ymajorgrids=true,
        xmajorgrids=true,
        grid style=dashed,
        unbounded coords=discard,
        filter discard warning=false,
        cycle list name=paulcolors6,
        legend style={at={(1.8,-3.1)}, anchor=north, draw=none, fill=none, legend columns=5},
    ]

    \nextgroupplot[title={hyperbolic}, ylabel={DG}]
    \addplot+[discard if not={solver}{sweep},discard if not={overlap}{1},mark repeat=20,mark size=2pt]
        table [x=steps, y=relative_residual, col sep=comma] {./numex/helmholtz_hyperbolic_dg_add.csv};
    \addplot+[discard if not={solver}{cg_sweep},discard if not={overlap}{1},mark repeat=20,mark size=2pt]
        table [x=steps, y=relative_residual, col sep=comma] {./numex/helmholtz_hyperbolic_dg_add.csv};
    \addplot+[discard if not={solver}{sweep},discard if not={overlap}{0},mark repeat=20,mark size=2pt]
        table [x=steps, y=relative_residual, col sep=comma] {./numex/helmholtz_hyperbolic_dg_add.csv};
    \addplot+[discard if not={solver}{cg_sweep},discard if not={overlap}{0},mark repeat=20,mark size=2pt]
        table [x=steps, y=relative_residual, col sep=comma] {./numex/helmholtz_hyperbolic_dg_add.csv};
    \legend{fixed point, PCG, fixed point (without overlap), PCG (without overlap)}
    \nextgroupplot[title={parabolic}]
    \addplot+[discard if not={solver}{sweep},discard if not={overlap}{1},mark repeat=20,mark size=2pt]
        table [x=steps, y=relative_residual, col sep=comma] {./numex/helmholtz_parabolic_dg_add.csv};
    \addplot+[discard if not={solver}{cg_sweep},discard if not={overlap}{1},mark repeat=20,mark size=2pt]
        table [x=steps, y=relative_residual, col sep=comma] {./numex/helmholtz_parabolic_dg_add.csv};
    \addplot+[discard if not={solver}{sweep},discard if not={overlap}{0},mark repeat=20,mark size=2pt]
        table [x=steps, y=relative_residual, col sep=comma] {./numex/helmholtz_parabolic_dg_add.csv};
    \addplot+[discard if not={solver}{cg_sweep},discard if not={overlap}{0},mark repeat=20,mark size=2pt]
        table [x=steps, y=relative_residual, col sep=comma] {./numex/helmholtz_parabolic_dg_add.csv};
    \nextgroupplot[title={elliptic}]
    \addplot+[discard if not={solver}{sweep},discard if not={overlap}{1},mark repeat=20,mark size=2pt]
        table [x=steps, y=relative_residual, col sep=comma] {./numex/helmholtz_elliptic_dg_add.csv};
    \addplot+[discard if not={solver}{cg_sweep},discard if not={overlap}{1},mark repeat=20,mark size=2pt]
        table [x=steps, y=relative_residual, col sep=comma] {./numex/helmholtz_elliptic_dg_add.csv};
    \addplot+[discard if not={solver}{sweep},discard if not={overlap}{0},mark repeat=20,mark size=2pt]
        table [x=steps, y=relative_residual, col sep=comma] {./numex/helmholtz_elliptic_dg_add.csv};
    \addplot+[discard if not={solver}{cg_sweep},discard if not={overlap}{0},mark repeat=20,mark size=2pt]
        table [x=steps, y=relative_residual, col sep=comma] {./numex/helmholtz_elliptic_dg_add.csv};
    
    \nextgroupplot[ylabel={TDG}]
    \addplot+[discard if not={solver}{sweep},discard if not={overlap}{1},mark repeat=20,mark size=2pt]
        table [x=steps, y=relative_residual, col sep=comma] {./numex/helmholtz_hyperbolic_tdg_add.csv};
    \addplot+[discard if not={solver}{cg_sweep},discard if not={overlap}{1},mark repeat=20,mark size=2pt]
        table [x=steps, y=relative_residual, col sep=comma] {./numex/helmholtz_hyperbolic_tdg_add.csv};
    \addplot+[discard if not={solver}{sweep},discard if not={overlap}{0},mark repeat=20,mark size=2pt]
        table [x=steps, y=relative_residual, col sep=comma] {./numex/helmholtz_hyperbolic_tdg_add.csv};
    \addplot+[discard if not={solver}{cg_sweep},discard if not={overlap}{0},mark repeat=20,mark size=2pt]
        table [x=steps, y=relative_residual, col sep=comma] {./numex/helmholtz_hyperbolic_tdg_add.csv};
    \nextgroupplot[]
    \addplot+[discard if not={solver}{sweep},discard if not={overlap}{1},mark repeat=20,mark size=2pt]
        table [x=steps, y=relative_residual, col sep=comma] {./numex/helmholtz_parabolic_tdg_add.csv};
    \addplot+[discard if not={solver}{cg_sweep},discard if not={overlap}{1},mark repeat=20,mark size=2pt]
        table [x=steps, y=relative_residual, col sep=comma] {./numex/helmholtz_parabolic_tdg_add.csv};
    \addplot+[discard if not={solver}{sweep},discard if not={overlap}{0},mark repeat=20,mark size=2pt]
        table [x=steps, y=relative_residual, col sep=comma] {./numex/helmholtz_parabolic_tdg_add.csv};
    \addplot+[discard if not={solver}{cg_sweep},discard if not={overlap}{0},mark repeat=20,mark size=2pt]
        table [x=steps, y=relative_residual, col sep=comma] {./numex/helmholtz_parabolic_tdg_add.csv};
    \nextgroupplot[]
    \addplot+[discard if not={solver}{sweep},discard if not={overlap}{1},mark repeat=20,mark size=2pt]
        table [x=steps, y=relative_residual, col sep=comma] {./numex/helmholtz_elliptic_tdg_add.csv};
    \addplot+[discard if not={solver}{cg_sweep},discard if not={overlap}{1},mark repeat=20,mark size=2pt]
        table [x=steps, y=relative_residual, col sep=comma] {./numex/helmholtz_elliptic_tdg_add.csv};
    \addplot+[discard if not={solver}{sweep},discard if not={overlap}{0},mark repeat=20,mark size=2pt]
        table [x=steps, y=relative_residual, col sep=comma] {./numex/helmholtz_elliptic_tdg_add.csv};
    \addplot+[discard if not={solver}{cg_sweep},discard if not={overlap}{0},mark repeat=20,mark size=2pt]
        table [x=steps, y=relative_residual, col sep=comma] {./numex/helmholtz_elliptic_tdg_add.csv};
    \nextgroupplot[ylabel={HDG}, xlabel={iterations}]
    \addplot+[discard if not={solver}{sweep},mark repeat=20,mark size=2pt]
        table [x=steps, y=relative_residual, col sep=comma] {./numex/helmholtz_hyperbolic_hdg_add.csv};
    \addplot+[discard if not={solver}{cg_sweep},mark repeat=20,mark size=2pt]
        table [x=steps, y=relative_residual, col sep=comma] {./numex/helmholtz_hyperbolic_hdg_add.csv};
    \nextgroupplot[xlabel={iterations}]
    \addplot+[discard if not={solver}{sweep},mark repeat=20,mark size=2pt]
        table [x=steps, y=relative_residual, col sep=comma] {./numex/helmholtz_parabolic_hdg_add.csv};
    \addplot+[discard if not={solver}{cg_sweep},mark repeat=20,mark size=2pt]
        table [x=steps, y=relative_residual, col sep=comma] {./numex/helmholtz_parabolic_hdg_add.csv};
    \nextgroupplot[xlabel={iterations}]
    \addplot+[discard if not={solver}{sweep},mark repeat=20,mark size=2pt]
        table [x=steps, y=relative_residual, col sep=comma] {./numex/helmholtz_elliptic_hdg_add.csv};
    \addplot+[discard if not={solver}{cg_sweep},mark repeat=20,mark size=2pt]
        table [x=steps, y=relative_residual, col sep=comma] {./numex/helmholtz_elliptic_hdg_add.csv};
    \end{groupplot} 
\end{tikzpicture}}
\vspace{-1.75em}
\caption{
    Fixed point and PCG with \textit{additive} preconditioners on the three trapping geometries (columns) for the three methods (rows).
}
\label{fig:helmholtz-trap3-preconditioner-add}
\end{figure}
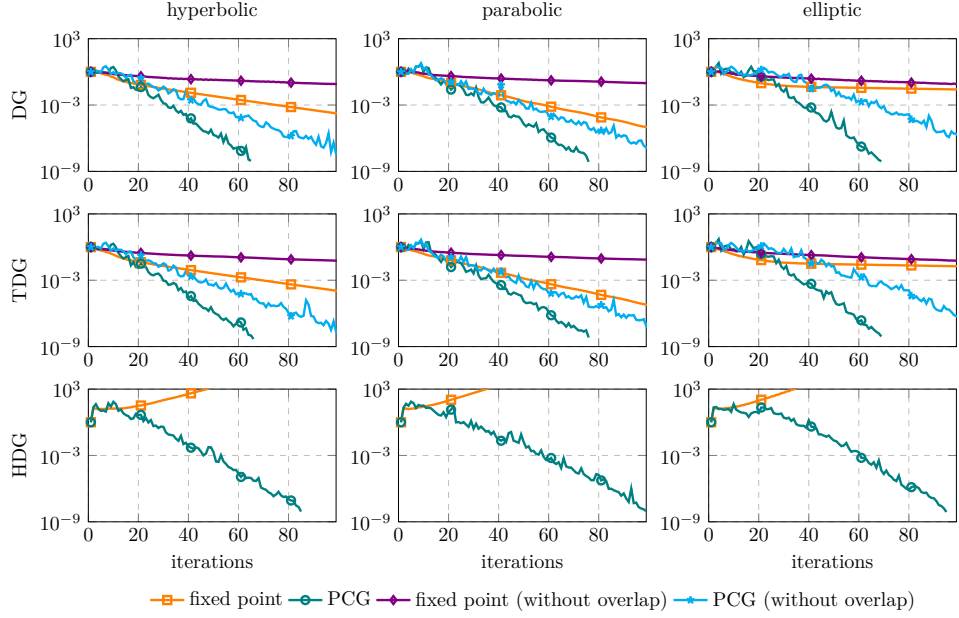
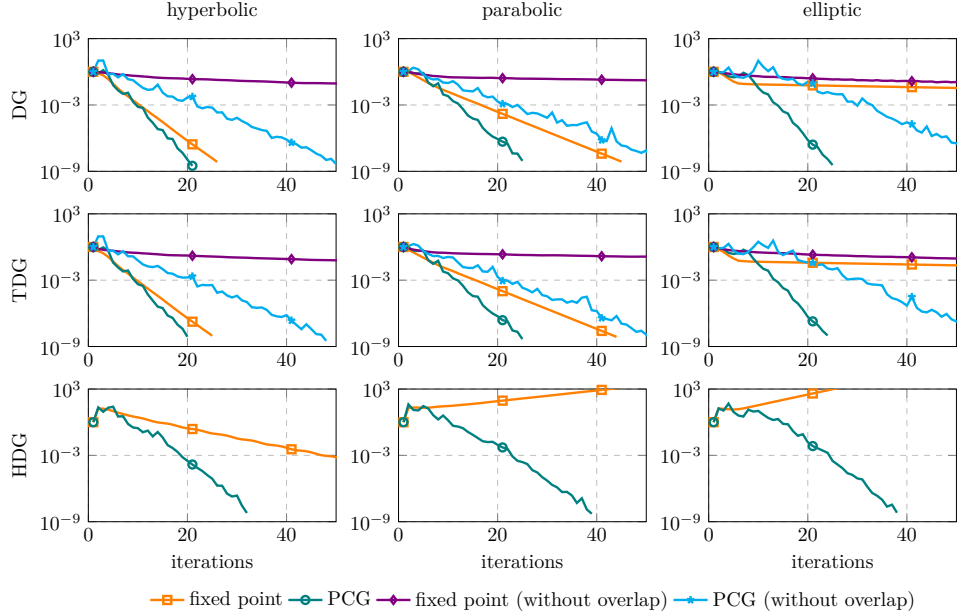
\begin{figure}[ht!]
\centering
\resizebox{\linewidth}{!}{
\begin{tikzpicture}
    \begin{groupplot}[
        group style={group size=3 by 3, horizontal sep=1.1cm, vertical sep=0.75cm},
        width=.47\linewidth,
        height=.31\linewidth,
        ymode=log,
        xmin=0,
        xmax=50,
        ymin=1e-9,
        ymax=1e3,
        ymajorgrids=true,
        xmajorgrids=true,
        grid style=dashed,
        unbounded coords=discard,
        filter discard warning=false,
        cycle list name=paulcolors6,
        legend style={at={(1.8,-3.1)}, anchor=north, draw=none, fill=none, legend columns=5},
    ]

    \nextgroupplot[title={hyperbolic}, ylabel={DG}]
    \addplot+[discard if not={solver}{sweep},discard if not={overlap}{1},mark repeat=20,mark size=2pt]
        table [x=steps, y=relative_residual, col sep=comma] {./numex/helmholtz_hyperbolic_dg_mult.csv};
    \addplot+[discard if not={solver}{cg_sweep},discard if not={overlap}{1},mark repeat=20,mark size=2pt]
        table [x=steps, y=relative_residual, col sep=comma] {./numex/helmholtz_hyperbolic_dg_mult.csv};
    \addplot+[discard if not={solver}{sweep},discard if not={overlap}{0},mark repeat=20,mark size=2pt]
        table [x=steps, y=relative_residual, col sep=comma] {./numex/helmholtz_hyperbolic_dg_mult.csv};
    \addplot+[discard if not={solver}{cg_sweep},discard if not={overlap}{0},mark repeat=20,mark size=2pt]
        table [x=steps, y=relative_residual, col sep=comma] {./numex/helmholtz_hyperbolic_dg_mult.csv};
    \legend{fixed point, PCG, fixed point (without overlap), PCG (without overlap)}
    \nextgroupplot[title={parabolic}]
    \addplot+[discard if not={solver}{sweep},discard if not={overlap}{1},mark repeat=20,mark size=2pt]
        table [x=steps, y=relative_residual, col sep=comma] {./numex/helmholtz_parabolic_dg_mult.csv};
    \addplot+[discard if not={solver}{cg_sweep},discard if not={overlap}{1},mark repeat=20,mark size=2pt]
        table [x=steps, y=relative_residual, col sep=comma] {./numex/helmholtz_parabolic_dg_mult.csv};
    \addplot+[discard if not={solver}{sweep},discard if not={overlap}{0},mark repeat=20,mark size=2pt]
        table [x=steps, y=relative_residual, col sep=comma] {./numex/helmholtz_parabolic_dg_mult.csv};
    \addplot+[discard if not={solver}{cg_sweep},discard if not={overlap}{0},mark repeat=20,mark size=2pt]
        table [x=steps, y=relative_residual, col sep=comma] {./numex/helmholtz_parabolic_dg_mult.csv};
    \nextgroupplot[title={elliptic}]
    \addplot+[discard if not={solver}{sweep},discard if not={overlap}{1},mark repeat=20,mark size=2pt]
        table [x=steps, y=relative_residual, col sep=comma] {./numex/helmholtz_elliptic_dg_mult.csv};
    \addplot+[discard if not={solver}{cg_sweep},discard if not={overlap}{1},mark repeat=20,mark size=2pt]
        table [x=steps, y=relative_residual, col sep=comma] {./numex/helmholtz_elliptic_dg_mult.csv};
    \addplot+[discard if not={solver}{sweep},discard if not={overlap}{0},mark repeat=20,mark size=2pt]
        table [x=steps, y=relative_residual, col sep=comma] {./numex/helmholtz_elliptic_dg_mult.csv};
    \addplot+[discard if not={solver}{cg_sweep},discard if not={overlap}{0},mark repeat=20,mark size=2pt]
        table [x=steps, y=relative_residual, col sep=comma] {./numex/helmholtz_elliptic_dg_mult.csv};
    
    \nextgroupplot[ylabel={TDG}]
    \addplot+[discard if not={solver}{sweep},discard if not={overlap}{1},mark repeat=20,mark size=2pt]
        table [x=steps, y=relative_residual, col sep=comma] {./numex/helmholtz_hyperbolic_tdg_mult.csv};
    \addplot+[discard if not={solver}{cg_sweep},discard if not={overlap}{1},mark repeat=20,mark size=2pt]
        table [x=steps, y=relative_residual, col sep=comma] {./numex/helmholtz_hyperbolic_tdg_mult.csv};
    \addplot+[discard if not={solver}{sweep},discard if not={overlap}{0},mark repeat=20,mark size=2pt]
        table [x=steps, y=relative_residual, col sep=comma] {./numex/helmholtz_hyperbolic_tdg_mult.csv};
    \addplot+[discard if not={solver}{cg_sweep},discard if not={overlap}{0},mark repeat=20,mark size=2pt]
        table [x=steps, y=relative_residual, col sep=comma] {./numex/helmholtz_hyperbolic_tdg_mult.csv};
    \nextgroupplot[]
    \addplot+[discard if not={solver}{sweep},discard if not={overlap}{1},mark repeat=20,mark size=2pt]
        table [x=steps, y=relative_residual, col sep=comma] {./numex/helmholtz_parabolic_tdg_mult.csv};
    \addplot+[discard if not={solver}{cg_sweep},discard if not={overlap}{1},mark repeat=20,mark size=2pt]
        table [x=steps, y=relative_residual, col sep=comma] {./numex/helmholtz_parabolic_tdg_mult.csv};
    \addplot+[discard if not={solver}{sweep},discard if not={overlap}{0},mark repeat=20,mark size=2pt]
        table [x=steps, y=relative_residual, col sep=comma] {./numex/helmholtz_parabolic_tdg_mult.csv};
    \addplot+[discard if not={solver}{cg_sweep},discard if not={overlap}{0},mark repeat=20,mark size=2pt]
        table [x=steps, y=relative_residual, col sep=comma] {./numex/helmholtz_parabolic_tdg_mult.csv};
    \nextgroupplot[]
    \addplot+[discard if not={solver}{sweep},discard if not={overlap}{1},mark repeat=20,mark size=2pt]
        table [x=steps, y=relative_residual, col sep=comma] {./numex/helmholtz_elliptic_tdg_mult.csv};
    \addplot+[discard if not={solver}{cg_sweep},discard if not={overlap}{1},mark repeat=20,mark size=2pt]
        table [x=steps, y=relative_residual, col sep=comma] {./numex/helmholtz_elliptic_tdg_mult.csv};
    \addplot+[discard if not={solver}{sweep},discard if not={overlap}{0},mark repeat=20,mark size=2pt]
        table [x=steps, y=relative_residual, col sep=comma] {./numex/helmholtz_elliptic_tdg_mult.csv};
    \addplot+[discard if not={solver}{cg_sweep},discard if not={overlap}{0},mark repeat=20,mark size=2pt]
        table [x=steps, y=relative_residual, col sep=comma] {./numex/helmholtz_elliptic_tdg_mult.csv};
    \nextgroupplot[ylabel={HDG}, xlabel={iterations}]
    \addplot+[discard if not={solver}{sweep},mark repeat=20,mark size=2pt]
        table [x=steps, y=relative_residual, col sep=comma] {./numex/helmholtz_hyperbolic_hdg_mult.csv};
    \addplot+[discard if not={solver}{cg_sweep},mark repeat=20,mark size=2pt]
        table [x=steps, y=relative_residual, col sep=comma] {./numex/helmholtz_hyperbolic_hdg_mult.csv};
    \nextgroupplot[xlabel={iterations}]
    \addplot+[discard if not={solver}{sweep},mark repeat=20,mark size=2pt]
        table [x=steps, y=relative_residual, col sep=comma] {./numex/helmholtz_parabolic_hdg_mult.csv};
    \addplot+[discard if not={solver}{cg_sweep},mark repeat=20,mark size=2pt]
        table [x=steps, y=relative_residual, col sep=comma] {./numex/helmholtz_parabolic_hdg_mult.csv};
    \nextgroupplot[xlabel={iterations}]
    \addplot+[discard if not={solver}{sweep},mark repeat=20,mark size=2pt]
        table [x=steps, y=relative_residual, col sep=comma] {./numex/helmholtz_elliptic_hdg_mult.csv};
    \addplot+[discard if not={solver}{cg_sweep},mark repeat=20,mark size=2pt]
        table [x=steps, y=relative_residual, col sep=comma] {./numex/helmholtz_elliptic_hdg_mult.csv};
    \end{groupplot}
\end{tikzpicture}}
\vspace{-1.5em}
\caption{
    Fixed point and PCG with \textit{multiplicative} preconditioners on the three trapping geometries (columns) for the three methods (rows).
}
\label{fig:helmholtz-trap3-preconditioner-mult}
\end{figure}

\Cref{fig:helmholtz-trap3-preconditioner-add} and \Cref{fig:helmholtz-trap3-preconditioner-mult} report the results for additive and multiplicative preconditioners respectively.
For each formulation, we compare fixed point iterative solver and PCG.
The partition of unity is only used in the additive version of the preconditioner with the fixed point solver (otherwise no convergence is obtained).

Across the trapping examples, the residuals show that the fixed point iterations are much more sensitive to trapped wave energy and that the Krylov-acceleration helps.
Overlap consistently improves the behavior of DG and TDG, while HDG remains competitive when used inside a Krylov iteration but seems reliable as a pure fixed point method.
The similar behavior across elliptic, parabolic, and hyperbolic trapping indicates that the trends are not tied to a single geometric configuration.

\subsection{Realistic submarine}

As a realistic three-dimensional scattering test we use the BeTSSi submarine benchmark, see \cite{SBGKMV_ICSAV_2003}.
The computation uses the CAD geometry scaled to a length of $62\,\mathrm{m}$ and water sound speed $1500\,\mathrm{m}/\mathrm{s}$.
We solve for the scattered field $u_s$ in a truncated exterior domain.
The incident plane wave $u_{\mathrm{inc}}$ is imposed through the sound-hard condition on the hull, $\partial_n u_s=-\partial_n u_{\mathrm{inc}}$.
The outer sphere uses a homogeneous first-order absorbing Robin condition for $u_s$.
The convergence tolerance is set to $10^{-8}$.
These runs use all 48 hardware threads available on the machine.

\begin{table}[ht!]
\centering
\footnotesize
\resizebox{\linewidth}{!}{%
\begin{filecontents*}{./numex/submarine_f100Hz_timings.csv}
case,method,solver,frequency_tag,overlap,pou,alpha,beta,frequency_hz,wavenumber,wavelength,order,elements_per_wavelength,ne,nv,ndofs,active_ndofs,ndomains,ncuts,min_domain_ndofs,avg_domain_ndofs,max_domain_ndofs,mat_nze,steps,residual,relative_residual,maxh,mesh_size,avg_domain_nnze,assembly_time,partition_time,doflist_time,preconditioner_time,setup_time,solve_time,avgmatvec_time,postprocess_time,tag,sound_speed,domain_radius,submarine_length,mesh_time,tol,maxiter,threads,cpu_count,time_to_solution,total_time
full,DG,cg_sweep,f100Hz,\cmark,\xmark,0.5,0.5,100.0,0.4188790204786391,14.999999999999998,5,2.0,150455,29079,8425480,8425480,100,16499,90832,102133.36,112280,2313392704,523,7.428404713616211e-09,7.428404713616211e-09,7.499999999999999,3990890.42624828,23133927.04,116.41224917198997,5.219500217004679,14.14261786796851,4836.753219658043,4856.116187373991,11303.447611747019,21.612710538713227,3.965978976339102e-06,full_dg,1500.0,50.0,62.0,9.330172647954896,1e-08,700,48,48,16301.763162752963,16302.895737671002
full,TDG,cg_sweep,f100Hz,\cmark,\xmark,0.5,0.5,100.0,0.4188790204786391,14.999999999999998,5,2.0,150455,29079,5416380,5416380,100,16499,58392,65657.16,72180,956044944,242,6.834693699985546e-09,6.834693699985546e-09,7.499999999999999,3990890.426248279,9560449.44,118.79640436096815,4.536601483006962,24.441029595967848,1926.4965307509992,1955.47472745704,2490.0787046630285,10.289581424227391,0.2654780619777739,full_tdg,1500.0,50.0,62.0,9.325198950013146,1e-08,700,48,48,4580.3698590809945,4581.577026292973
full,HDG,cg_sweep,f100Hz,\cmark,\xmark,1.0,1.0,100.0,0.4188790204786391,14.999999999999998,5,2.0,150455,29079,21282490,9246090,100,16499,93300,97410.6,102210,1902296700,314,8.473491088975698e-09,8.473491088975698e-09,7.499999999999999,3990890.42624828,19022967.0,244.93971583002713,4.413411077985074,43.03485802101204,3584.461009390012,3631.9099290669546,4623.1088230099995,14.723276506401271,11.38200477999635,full_hdg,1500.0,50.0,62.0,9.477296018972993,1e-08,700,48,48,8535.743673608988,8537.987420881982
\end{filecontents*}
\pgfplotstabletypeset[
    col sep=comma, fixed, precision=2,
    columns={
    method,overlap,pou,active_ndofs,assembly_time,setup_time,
    solve_time,postprocess_time,total_time,steps
    },
    columns/method/.style={column name={method}, column type=c, string type},
    columns/overlap/.style={column name={overlap}, column type=c, string type},
    columns/pou/.style={column name={PoU}, column type=c, string type},
    columns/active_ndofs/.style={column name={ndofs}, column type=c},
    columns/assembly_time/.style={column name={assembly}, column type=c},
    columns/setup_time/.style={column name={factorization}, column type=c},
    columns/solve_time/.style={column name={solve}, column type=c},
    columns/postprocess_time/.style={column name={postprocess}, column type=c},
    columns/total_time/.style={column name={total}, column type=c},
    columns/steps/.style={column name={iterations}, column type=c},
    every head row/.style={
        before row=\toprule,
        after row=\midrule,
    },
    every odd row/.style={
        before row={\rowcolor{gray!15}}
    },
    every last row/.style={after row=\bottomrule},
]{./numex/submarine_f100Hz_timings.csv}
}
\vspace{.5em}
\caption{
Timings in seconds for the BeTSSi submarine benchmark at \(f=100\,\mathrm{Hz}\), \(p=5\), two elements per wavelength, and \(J=100\) subdomains.
}
\label{tab:submarine-timings}
\vspace{-1em}
\end{table}

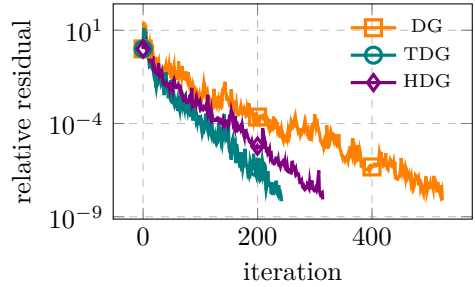
\begin{figure}[ht!]
\begin{tikzpicture}
\begin{semilogyaxis}[
    width=.5\linewidth,
    height=.35\linewidth,
    xlabel={iteration},
    ylabel={relative residual},
    ymajorgrids=true,
    xmajorgrids=true,
    grid style=dashed,
    legend style={draw=none, fill=none, font=\scriptsize, at={(0.98,0.98)}, anchor=north east},
    cycle list name=paulcolors6,
]
\addplot+[mark repeat=200]
    table [x=iteration, y=relative_residual, col sep=comma]
    {./numex/submarine_dg_f100Hz_history.csv};
\addlegendentry{DG}
\addplot+[mark repeat=200]
    table [x=iteration, y=relative_residual, col sep=comma]
    {./numex/submarine_tdg_f100Hz_history.csv};
\addlegendentry{TDG}
\addplot+[mark repeat=200]
    table [x=iteration, y=relative_residual, col sep=comma]
    {./numex/submarine_hdg_f100Hz_history.csv};
\addlegendentry{HDG}
\end{semilogyaxis}
\end{tikzpicture}
\caption{
Convergence history for the BeTSSi submarine benchmark.
}
\label{fig:submarine-history}
\end{figure}

For the first run, the acoustic frequency is chosen as $100\,\mathrm{Hz}$, corresponding to $\omega\approx0.418\,\mathrm{m}^{-1}$.
We use polynomial order \(p=5\), two elements per wavelength, and \(J=100\) subdomains.
In \Cref{tab:submarine-timings} we report timings for the different methods and show the convergence history in \Cref{fig:submarine-history}, using a multiplicative preconditioner with overlap and the CG solver.
This proof-of-concept implementation performs the subdomain factorizations and multiplicative sweeps sequentially, so the timings should not be read as optimized parallel DD timings.
Notably, the TDG assembly time is comparable to DG even though it includes Trefftz-space setup via element-wise SVD, which indicates that this local construction parallelizes well.
As before, the timings identify the linear solve and the preconditioner setup and application as the dominant parts of the computation.
While the number of non-zero entries of HDG will be smaller than DG, the globally coupled degrees of freedom (ndof in the table) are larger than for DG and TDG, recall \Cref{tab:helmholtz-simplicial-counts}.
Again, we only consider the degrees of freedom that remain after condensation to the trace system.
The residual history, shown in \Cref{fig:submarine-history}, is less regular than in the plane-wave tests.

The submarine benchmark mainly demonstrates how the solver behavior changes once the smooth manufactured setting is replaced by a realistic exterior scattering geometry.
This makes reductions in local problem size and iteration cost especially important, and suggests that the conclusions from the controlled tests remain relevant at the scale of a practical scattering computation.

\begin{table}[ht!]
\centering
\footnotesize
\resizebox{\linewidth}{!}{%
\begin{filecontents*}{./numex/submarine_f200Hz_timings.csv}
case,method,solver,frequency_tag,overlap,pou,alpha,beta,frequency_hz,wavenumber,wavelength,order,elements_per_wavelength,ne,nv,ndofs,active_ndofs,ndomains,ncuts,min_domain_ndofs,avg_domain_ndofs,max_domain_ndofs,mat_nze,steps,residual,relative_residual,maxh,mesh_size,avg_domain_nnze,assembly_time,partition_time,doflist_time,preconditioner_time,setup_time,solve_time,avgmatvec_time,postprocess_time,tag,sound_speed,domain_radius,submarine_length,mesh_time,tol,maxiter,threads,cpu_count,time_to_solution,total_time
full,TDG,cg_sweep,f200Hz,\cmark,\xmark,0.5,0.5,200.0,0.8377580409572782,7.499999999999999,6,2.0,183051,35586,8969499,8969499,400,34652,25872,30670.325,34496,2153296033,205,8.13077718694508e-09,8.13077718694508e-09,3.7499999999999996,2223783.6362431385,5383240.0825,408.0116038520355,5.8941764569608495,33.631031602970324,12807.43705268699,12846.963720929052,3846.3115293129813,18.76249526494137,0.6535345729789697,full_tdg,1500.0,50.0,62.0,11.695197905995883,1e-08,700,48,48,17126.73305430502,17128.39259464998
\end{filecontents*}
\pgfplotstabletypeset[
    col sep=comma, fixed, precision=2,
    columns={
    method,overlap,pou,active_ndofs,assembly_time,setup_time,
    solve_time,postprocess_time,total_time,steps
    },
    columns/method/.style={column name={method}, column type=c, string type},
    columns/overlap/.style={column name={overlap}, column type=c, string type},
    columns/pou/.style={column name={PoU}, column type=c, string type},
    columns/active_ndofs/.style={column name={ndofs}, column type=c},
    columns/assembly_time/.style={column name={assembly}, column type=c},
    columns/setup_time/.style={column name={factorization}, column type=c},
    columns/solve_time/.style={column name={solve}, column type=c},
    columns/postprocess_time/.style={column name={postprocess}, column type=c},
    columns/total_time/.style={column name={total}, column type=c},
    columns/steps/.style={column name={iterations}, column type=c},
    every head row/.style={
        before row=\toprule,
        after row=\midrule,
    },
    every odd row/.style={
        before row={\rowcolor{gray!15}}
    },
    every last row/.style={after row=\bottomrule},
]{./numex/submarine_f200Hz_timings.csv}
}
\vspace{.5em}
\caption{
Timings in seconds for the BeTSSi submarine benchmark at \(f=200\,\mathrm{Hz}\), \(p=6\), two elements per wavelength, and \(J=400\) subdomains.
}
\label{tab:submarine-timings2}
\vspace{-1em}
\end{table}

\begin{figure}[ht!]
\centering
\adjustbox{valign=c}{%
\maybeimage[width=.49\linewidth,clip,trim=17.6cm 9cm 9.5cm 9cm]{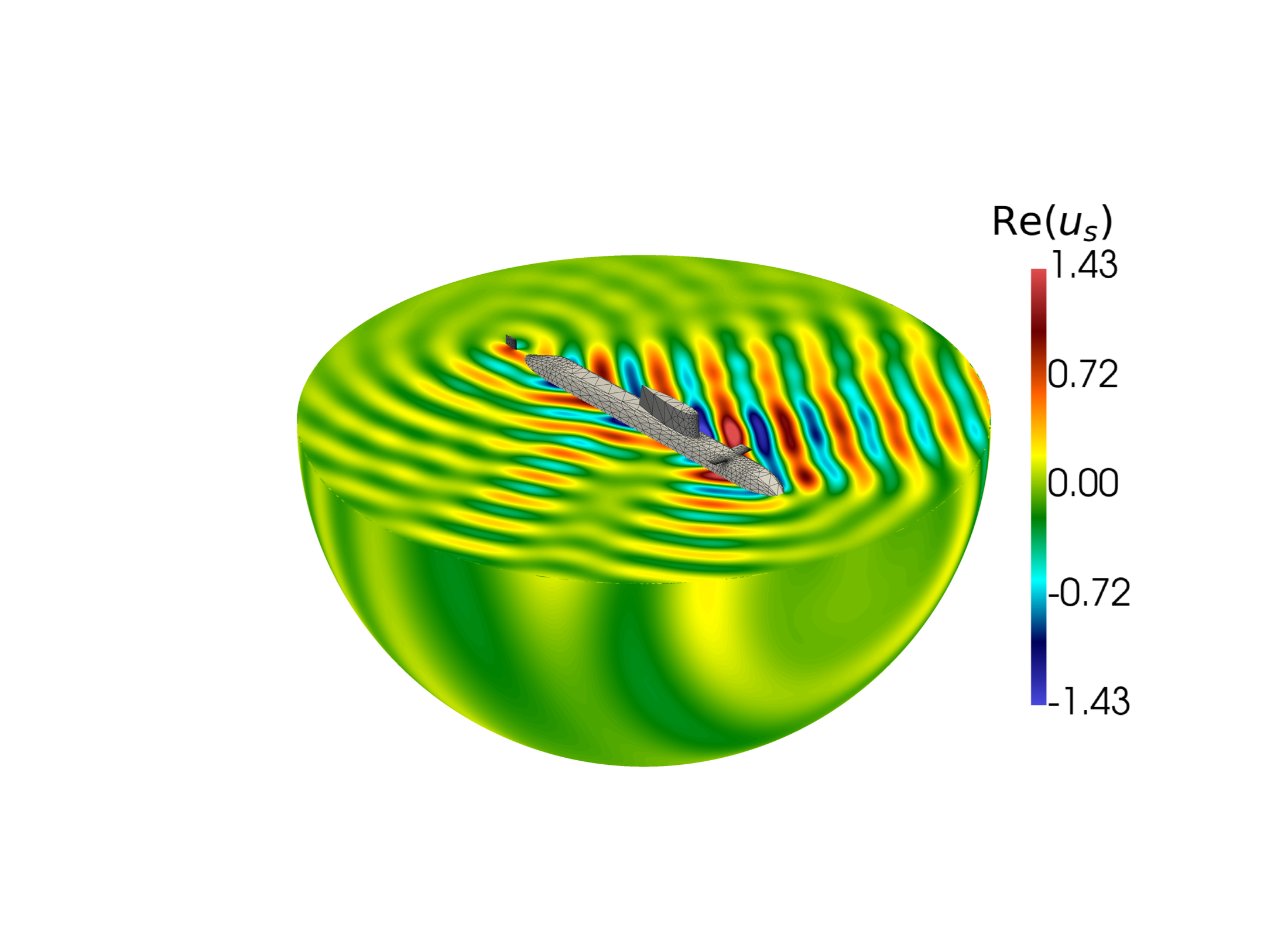}%
}
\hfill%
\adjustbox{valign=c}{%
\maybeimage[width=.49\linewidth,clip,trim=14cm 1cm 8cm 1cm]{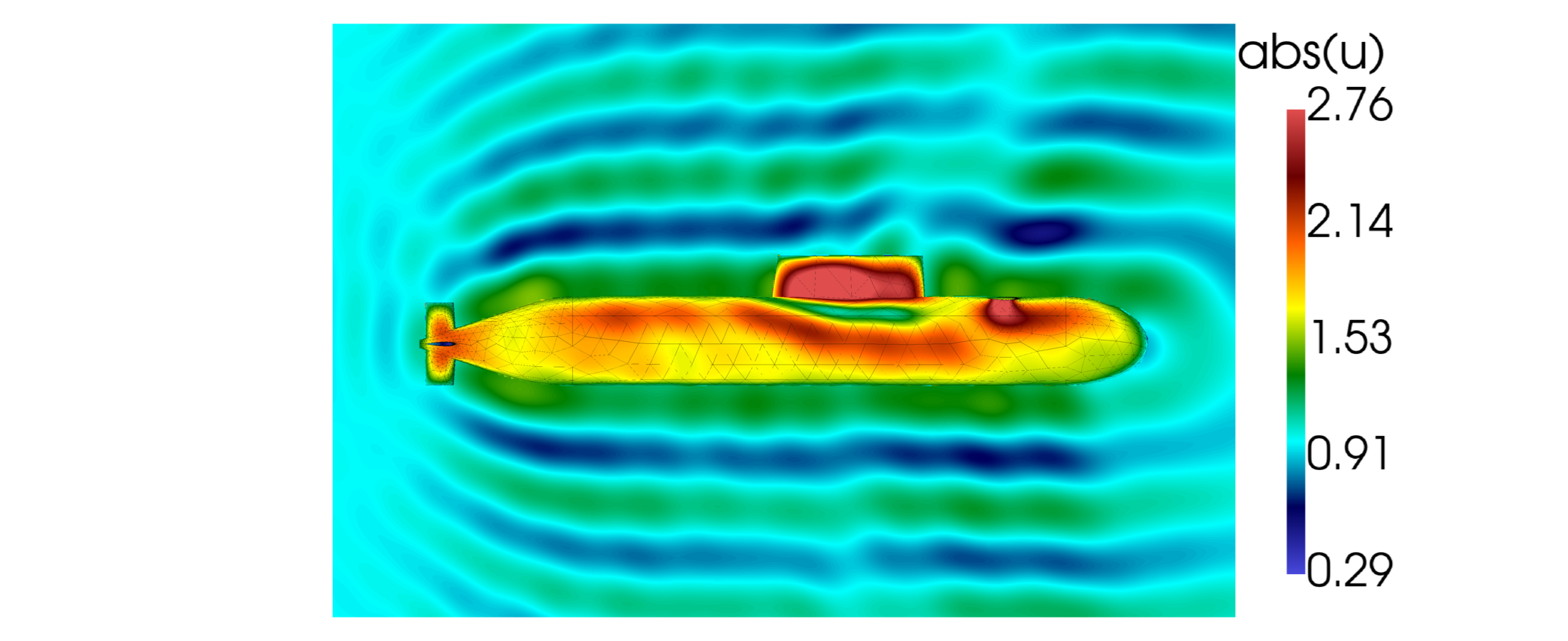}%
}
\caption{
Visualization of the BeTSSi submarine benchmark with TDG, $p=6$, and $f=200\,\mathrm{Hz}$.
}
\label{fig:submarine-visualization}
\end{figure}

In \Cref{tab:submarine-timings2} we show timings for the same benchmark for $200\,\mathrm{Hz}$, corresponding to $\omega\approx0.838\,\mathrm{m}^{-1}$, with \(p=6\) and \(J=400\) subdomains.
Compared with the \(100\,\mathrm{Hz}\) TDG run, the larger run mainly increases the factorization time, partly because the subdomain factorizations are performed sequentially.
The iteration count is lower than for the \(100\,\mathrm{Hz}\) TDG run, but the comparison is not a pure frequency study because the polynomial order, mesh, and number of subdomains also change.
For this choice of parameters we could only run the TDG method, as DG and HDG exceeded the available memory of 512 GB.
In \Cref{fig:submarine-visualization} we show a visualization of the TDG solution at \(f=200\,\mathrm{Hz}\).
The overview plots $\Re u_s$ to isolate the scattered response of the hull, while the close-up shows $|u_s+u_{\mathrm{inc}}|$, the total-pressure amplitude relevant near the obstacle.

\section*{Acknowledgments}
\sloppy{
This research was funded in part by the Austrian Science Fund (FWF) \href{https://doi.org/10.55776/ESP4389824}{10.55776/ESP4389824}.
For open access purposes, the author has applied a CC BY public copyright license to any author-accepted manuscript version arising from this submission.}

\bibliographystyle{amsplain}
\bibliography{bib}

\end{document}